\documentclass{tac}
\usepackage{xy}
\usepackage{MnSymbol}
\usepackage{amsfonts}

\usepackage{tikz}

\xyoption{all}

\def\mathcaldef#1{\expandafter\def\csname#1\endcsname{{\cal#1}}}

\def\"{``}
\def\q{\quad}
\def\qq{\quad\quad}
\def\qv{\qq ;\qq}

\def\iso{\,\cong\,}

\def\adj{\dashv}
\def\-{\text{-}}
\def\op{^{\rm op}}

\def\ov{\overline}

\def\tm{\times}

\def\al{\alpha}

\def\t{_\blacktriangleright}

\def\dSpan{\mathbb{S}\rm{pan}\,}
\def\dCat{\mathbb{C}\rm{at}}
\def\dSq{\mathbb{S}\rm{q}\,}
\def\dPb{\mathbb{P}\rm{b}\,}
\def\dTot{\mathbb{T}\rm{ot}\,}
\def\dBij{\mathbb{B}\rm{ij}\,}
\def\dEspan{\mathbb{E}\rm{span}}
\def\dSec{\mathbb{S}\rm{ec}}
\def\dA{\mathbb A}
\def\dM{\mathbb M}
\def\dN{\mathbb N}
\def\dB{\mathbb B}
\def\dC{\mathbb C}

\def\dH{\mathbb H\,}

\def\dP{\mathbb P\,}

\def\dG{\mathbb G\,}

\def\dR{\mathbb R\,}

\mathcaldef{A}
\mathcaldef{B}
\mathcaldef{C}
\mathcaldef{D}
\mathcaldef{S}
\mathcaldef{DP}
\mathcaldef{DG}
\mathcaldef{E}
\mathcaldef{L}

\mathbfdef{Set}
\mathbfdef{Set_f}
\mathbfdef{Bij}
\mathbfdef{Tot}
\mathbfdef{Mon}
\mathbfdef{cMon}
\mathbfdef{MonCat}
\mathbfdef{sMonCat}
\mathbfdef{Span}
\mathbfdef{Espan}
\mathbfdef{Cat}
\mathbfdef{sCat}
\mathbfdef{dCat}
\mathbfdef{sdCat}
\mathbfdef{Prof}
\mathrmdef{id}
\mathbfdef{sMlt}
\mathbfdef{Mlt}
\mathbfdef{Pb}
\mathbfdef{Sp}
\mathbfdef{colSp}
\mathbfdef{Sec}

\def\obj{{\rm obj\,}}
\def\pb{{\rm pb}}
\def\Fam{{\rm Fam\,}}
\def\Famf{{\rm Fam}_f\,}
\def\core{{\rm core\,}}
\def\cart{{\rm cart}}

\newtheorem{prop}{Proposition}

\let\pf\proof
\let\epf\endproof
\def\eq{\begin{equation}}
\def\eeq{\end{equation}}

\author{Claudio Pisani}
\address{via Saluzzo 67,\\ 10125 Torino, Italy.}
\title{Unbiased Multicategory Theory}
\copyrightyear{2024}
\eaddress{pisclau@yahoo.it}

\begin{document}

\maketitle

\begin{abstract}

We present an unbiased theory of symmetric multicategories, where sequences are replaced by families. 
To be effective, this approach requires an explicit consideration of indexing and reindexing of objects and arrows,
handled by the double category $\dPb$ of pullback squares in finite sets: 
a symmetric multicategory is a sum preserving discrete fibration of double categories $M: \dM\to \dPb$.
If the \"loose" part of $M$ is an opfibration we get unbiased symmetric monoidal categories.

The definition can be usefully generalized by replacing $\dPb$ with another double prop $\dP$,
as an indexing base, giving $\dP$-multicategories.
For instance, we can remove the finiteness condition to obtain infinitary symmetric multicategories, 
or enhance $\dPb$ by totally ordering the fibers of its loose arrows to obtain plain multicategories.

We show how several concepts and properties find a natural setting in this framework.
We also consider cartesian multicategories as algebras for a monad $(-)^\cart$ on $\sMlt$,
where the loose arrows of $\dM^\cart$ are \"spans" of a tight and a loose arrow in $\dM$.

\end{abstract}


\section{Introduction}
\label{intro}

We present an unbiased notion of multicategory,
by exploiting the language of double categories and of fibrations.
While an immediate advantage of such a notion is to remove the awkwardness which
inevitably derive from a skeletal approach, we hope to convince the reader that 
it also renders neater various related concepts and properties, opening up new perspectives. 

Here, by \"unbiased" we refer to the idea that, in defining a notion which involves families and operations on them, 
one should not impose any particular structure (such as an order) on the indexing sets of the families, 
unless such a structure is really part of the notion itself. 
These spurious structures come up usually when one uses a particular skeleton of the category of (finite) sets
(such as that formed by the ordinals) as an indexing for the families.
So, for instance, the idea of a commutative monoid $M$ should be that one can add the elements
of a family $I\to M$, regardless of any order on the set $I$ and, by reindexing the family through a bijection $J\to I$,
the addition of the new family $J\to M$ should give the same result.
(See Section \ref{Cm} for a detailed account of unbiased commutative monoids.)
While the price to pay is that one has to take in account the indexing and reindexing explicitly,
the unbiased approach turns out ultimately to have a simplifying role.

The same term is often also referred to the fact that in studying algebraic structures 
one should treat all the operations on the same footing, as in Lawvere theories.
For instance, in the theory of monoids it is often rather unnatural to insist on the generating set 
of binary and nullary operations, as it is usually done.
The two aspects are of course not unrelated: it would be unnatural to restrict to the families indexed by sets of given cardinalities.
Thus, the concepts defined here are actually unbiased also in this sense.

In fact, in the classical notion of multicategory one chooses to work  
with a skeleton of $\Set_f$ as an indexing for the families of objects and arrows.
Though this may superficially appear as a simplifying step, already in the definition itself 
composition and associativity are expressed in a rather intricate way.
Things get even worse for symmetric multicategories, where one needs to express also 
the compatibility of composition with permutations of objects and arrows, 
and where it is intuitively clear that ordering of objects or arrows should play no role.


\subsection{The role of double categories}
\label{Mi}\qq 
The main idea is to consider explicitly the indexing and reindexing of objects and arrows 
in a symmetric multicategory, 
via a double functor $\dM \to \dPb$ to the double category of pullback squares in finite sets.
The bureaucracy of composition and of compatibility with permutations (reindexing)
is then automatically handled by the double categoral setting.
Given a symmetric multicategory $M$ one can in fact use
\begin{itemize}
\item
the objects of $\dPb$, that are finite sets, to index the families of objects of $M$;
\item
the tight arrows of $\dPb$, that are mappings, to reindex the families of objects;
\item
the loose arrows of $\dPb$, that are mappings, to index the families of arrows of $M$;
\item
the cells of $\dPb$, that are pullback squares, to reindex the families of arrows.
\end{itemize}
We so get an obvious functor $\dM\to \dPb$, where $\dM$ is the double category such that
\begin{itemize}
\item
$\dM_0$ is the category of families of objects of $M$ with their reindexing;
\item
$\dM_1$ is the category of families of arrows of $M$ with their reindexing;
\item
loose composition is given by composition in $M$.
\end{itemize}

Thus, we are led to {\em define} a symmetric multicategory as discrete fibration of double categories 
$M: \dM\to \dPb$ (meaning that both the components 
$M_0 : \dM_0 \to \Set_f$ and $M_1: \dM_1\to \dPb_1$ of $M$ are discrete fibrations).
To ensure that the objects and the loose arrows of $\dM$ are really families of single \"objects" and \"arrows",
we need also to assume that the categories $\dM_0$ and $\dM_1$ have finite sums and that they
preserved both by the source and target functors $s,t : \dM_1 \to \dM_0$ and by 
the components $M_0$ and $M_1$ of $M$.

We have so fully captured the notion of unbiased symmetric multicategory as a sum preserving
discrete fibration $M: \dM\to \dPb$.
A morphism $M\to M'$ is simply a sum preserving double functor $\dM\to\dM'$ over $\dPb$.


\subsection{Changing the base}
\label{Bas}\qq
This definition of unbiased symmetric multicategory can be usefully generalized by 
replacing $\dPb$ with another \"indexing base" $\dP$.

On the one hand, we can simply remove the finiteness condition 
(that is, replace $\Set_f$ with $\Set$) to obtain infinitary symmetric multicategories or monoidal categories.
On the other hand, we can use any double prop $\dP$ as the indexing base in place of $\dPb$.

Here, by a \"double prop" we mean a double category $\dP$ with finite sums such that $\dP_0 = \Set_f$, 
or $\dP_0 = \Set_f/S$ for colored double props (see Section \ref{Dc}).
Colored double props form a category $\DP$ whose arrows are sum preserving discrete fibrations.
As we have just sketched, the category $\sMlt$ of (unbiased) symmetric multicategories 
can be defined as the slice of $\DP$ over the double prop $\dPb$:
\eq 
\sMlt = \DP / \dPb 
\eeq

By replacing $\dPb$ with any double prop $\dP$ we get the category of $\dP$-multicategories:
\eq 
\dP\-\Mlt = \DP / \dP 
\eeq 
so that $\dPb$-multicategories are symmetric multicategories.

For instance, we can consider the double prop $\dTot$, whose loose arrows are mappings
endowed with a total order on each of its fibers, cells are again pullbacks 
(along which order are transported) and loose composition is given by juxtaposing orders.
In this case we get the category of (unbiased) plain (or planar) multicategories:
\eq
\dTot\-\Mlt = \Mlt
\eeq 

Furthermore, the forgetful functor $\dTot \to \dPb$, is itself a (one-object) symmetric multicategory $T \in \sMlt$ 
(classically known as the \"associative operad").
Thus, by a standard argument on slice categories,
\eq 
\sMlt / T = \Mlt 
\eeq

Similarly, if $\dBij$ is the double subcategory of $\dPb$ formed by those cells whose
loose sides are bijections, we get categories (since we have only unary arrows):
\eq \label{bas1}
\dBij\-\Mlt = \Cat
\eeq
Again, the inclusion $\dBij \to \dPb$ is itself a symmetric multicategory $U \in \sMlt$ (with one object and one arrow), 
so that we have the standard equivalence 
\eq
\sMlt / U = \Cat 
\eeq

More generally, the domain $\dM$ of any $\dP$-multicategory $M:\dM\to\dP$, 
is itself a double prop and we have the equivalence
\eq \label{sl} 
\dM\-\Mlt = (\dP\-\Mlt) / M 
\eeq
 
Consider for instance a category $\C$ as a $\dBij$-multicategory as in (\ref{bas1}): 
\eq \label{bas2}
\C:\dC\to\dBij
\eeq
(so that the loose part of $\dC$ is formed by families of arrows in $\C$).
Then the $\dC$-multicategories are the categories over $\C$:
\eq
\dC\-\Mlt = \Cat/\C
\eeq
 

\subsection{Unbiased monoidal categories}
\label{Unmon}\qq
Let $M: \dM\to \dPb$ be a symmetric multicategory.
If $M_l: \dM_l\to \Set_f$, the loose part of $M$, is an opfibration and its opcartesian arrows 
are stable with respect to reindexing, we get unbiased symmetric monoidal categories:
the opcartesian arrows give the tensor product in its universal or representable form, as in \cite{hermida}.
If $M_l: \dM_l\to \Set_f$ is a {\em discrete} opfibration we get unbiased commutative monoids
(see Section \ref{Cm}).
Similarly, for any double prop $\dP$, we define the category of $\dP$-monoidal categories
and the category of $\dP$-algebras (or $\dP$-monoids). 

Let us see some instances:
\begin{itemize}
\item
$\dTot$-monoidal categories are the (unbiased) monoidal categories, 
while $\dTot$-algebras are the (unbiased) monoids.
\item
If $\dC \in \DP$ corresponds to the category $\C$ as in (\ref{bas2}) above,
$\dC$-monoidal categories (respectively, $\dC$-algebras) are of course the opfibrations 
(respectively, the discrete opfibrations) over $\C$.
\item
If $\dM \in \DP$ is the double prop corresponding to a symmetric multicategory $M:\dM \to \dPb$, 
the category of $\dM$-algebras is the fibered form of the usual category of $M$-algebras,
that is, of functors from $M$ to the symmetric multicategory of sets and several variable mappings.
The same holds for a plain multicategory $M:\dM \to \dTot$.
\item
If $M:\dM \to \dP$ is a $\dP$-algebra, the $\dM$-algebras are simply $\dP$-algebras over $\dM$.
For instance, if $M$ is a (commutative) monoid, its algebras are (commutative) monoids over it.
\end{itemize}

Summarizing, any double prop $\dP \in \DP$ can be seen as a \"theory" and any morphism  
$F:\dP' \to \dP$ (that is, any $\dP$-multicategory) as a \"lax model" for $\dP$. 
So, the domain $\dP'$ of a lax model is itself a theory, 
and its models are also models for $\dP$ (via composition with $F$).
Furthermore, a model for $\dP$ can be \"representable" (giving $\dP$-monoidal categories)
or \"strict" (giving $\dP$-algebras).

In this respect, the prop for commutative monoids, (a skeleton of) $\Set_f$, 
is a trace of the double prop $\dPb$ for symmetric multicategories. 
Commutative monoids arise as the strict models.
Similarly, $\core(\Set_f)$, the prop for sets, is a trace of $\dBij$, the double prop for categories.
Sets arise as the strict models.
As for plain monoids, classically one considers the pro (rather than the prop)
of finite ordinals (a skeleton of the category of total orders). 
But they can also be obtained as the algebras for the (non-skeletal) prop whose
objects are finite sets, and whose arrows are mappings with a total order on each fiber. 
This is a trace of $\dTot$, the double prop for (plain) multicategories.
Again, monoids arise as the strict models.


\subsection{Outline}
\label{Out}\qq
In Section \ref{SET}, we recall and set up some technical tools to be used in our framework.
In particular, we consider the family construction, which provides the free sum completion of a category.

The category of finite sets and bijections, $\core(\Set_f)$, is the base category for expressing 
symmetry and transport of structures, as in Joyal theory of species.
Its free sum completion $\Fam(\core(\Set_f))$ is the category of $\Pb$ of mappings and pullback squares.
What is gained is that $\Pb$ is naturally the arrow part of a double category structure $\dPb$, 
which allows us to exploit the double category machinery.
In particular, we get the unbiased form of a \"colored symmetric sequence" (or colored species) as 
 a double graph with sums and with a sum preserving discrete fibration over $\dPb$.

In Section \ref{USM}, we can then consider unbiased symmetric multicategories as multiplicative species,
that is monoids in the category of colored symmetric sequences (with respect to the monoidal structure 
inherited from the double category structure on the double graph $\dPb$). 
Another point of view is to consider unbiased symmetric multicategories as double props $\dP$ 
with an \"indexing structure" (sum preserving discrete fibration) $M:\dP \to \dPb$.
Intuitively, this is the unbiased version of the presentation of multicategories in terms of the free props
on them given in \cite{leinster}. The rather {\em ad hoc} characterization of the latter given there, is here 
obtained by the \"indexing structure" $M:\dP \to \dPb$, which assures that maps over any mapping are indeed
families (sums) of maps with a single object as codomain.

In Section \ref{CART}, we define a monad $(-)^\cart$ on $\sMlt$, whose algebras are (unbiased)
cartesian multicategories.
That we so indeed capture the right notion of cartesian multicategory
is confirmed by the fact that one can prove therein the equivalence between 
tensor products (representability), universal products and algebraic products
(Proposition \ref{maincart}).


\subsection{Related work}
Most of the ideas presented here originate from the previous work of the author \cite{pisani2}
and \cite{pisani3}. In particular, the fibrational and double categorical approach to symmetric multicategories 
had been developed already in \cite{pisani2}, albeit in a rather rough form.
In \cite{pisani3}, the same framework is considered under a different perspective,
namely by exploiting lax double functors in place of fibrations of double categories.

Indeed (see \cite{lambert}), there is a correspondence between discrete fibrations of double categories $\dP' \to \dP$
and {\em lax} double functors $\dP \to \dSpan(\Set)$. 
Furthermore, if $\dP$ and $\dP'$ have sums, sum preserving discrete fibrations $\dP' \to \dP$ correspond to 
product preserving lax functors $\dP\op \to \dSpan(\Set)$.
In \cite{patterson}, finite product double categories are considered as theories and finite product-preserving lax functors
out of them as their models (with $\dSpan(\Set)$ as the default receiving semantics).
Thus, the above sketched idea of double props as theories may be seen as exhibiting them (actually, their dual) 
as instances of theories in the sense of \cite{patterson}.
What is remarkable is that, while in that paper most of the instances of theories are presented via generators and relations,
our instances are exhibited explicitly as pretty natural double categories and their models turn out to be 
important kinds of unbiased multicategory-like structures. 
Furthermore, the fibrational approach adopted here has the advantage of enlightening the aspects 
related to the change of base issues, as sketched in Section \ref{Bas}. 

An unbiased notion of symmetric multicategory, where sequences are replaced by families,
is already present in the literature, though the definition seems rather unwieldy;
see for instance \cite{beilinson} (\"pseudo-tensor categories"),  \cite{leinster} (\"fat symmetric multicategories") and \cite{lurie}).

On the other hand, there is the idea of associating, to a (symmetric) multicategory $M$, 
the ordinary category of finite families (or sequences) of arrows in $M$.
While we already mentioned the related free prop approach to multicategories in \cite{leinster},
the idea itself dates back at least to \cite{may} (see also \cite{lurie}), with the so called  
\"categories of operators". 
These are defined as categories over {\em pointed} finite sets (rather than over finite sets) and it is then possible to 
characterize symmetric multicategories with rather complicated conditions on these functors.
Again, in the author's opinion, this approach seems to lack clarity and flexibility. 

We believe that, to obtain a really useful notion of unbiased multicategory, one has also to explicitly consider the reindexing 
as part of the structure, exploiting the language of double categories and of fibrations.


\subsection{Aknowledgements}
\label{Ak}\qq 
I wish to express my gratitude to Nathanael Arkor, with whom I have had several fruitful exchanges of ideas.
His constant encouragement and his helpful remarks and suggestions 
have contributed greatly to improving this work and to giving it its present form.


\section{The setting}
\label{SET}

In the present section we recall and set up some technical tools to be used in our framework,
in particular concerning the family construction and double categories. 
We are then in a position to consider species in the sense of \cite{joyal} as double graphs $\dG$
with an \"indexing structure" $\dG\to\dPb$.


\subsection{The family construction}
\label{Fam}\qq 
The family construction (introduced in \cite{benabou2}) plays a prominent role in the present context 
and it is of overall importance in category theory.
We here briefly review some of its properties and of its instances,
most of which are probably well known. 
Maybe this is not the case for Proposition \ref{fam4}.

Given a category $\C$, we can construct another category $\Fam\C$
as follows:
\begin{itemize}
\item
an object of $\Fam\C$ is a family of objects of $\C$, that is a set $I$ along with a mapping $A:I\to\obj\C$, 
often denoted $A_i \, (i\in I)$;
\item
an arrow $A_i \, (i\in I) \to B_j \, (j\in J)$ is a mapping $f:I\to J$ along with 
a family of arrows $\al_i : A_i \to B_{fi} \, (i\in I)$ in $\C$.
\end{itemize}

By the construction itself, there is an obvious functor $\Fam\C \to \Set$,
which is a split fibration (the \"family fibration" of $\C$).
The cartesian arrows are those for which all the $\al_i$ are isomorphisms,
and a splitting is given by the arrows $\id_{B_i} : B_{fi} \to B_{fi} \, (i\in I)$.
It is a bifibration if and only if $\C$ has sums.

We denote by $\Famf\C$ the full subcategory whose objects are indexed by finite sets. 

It is well known that $\Fam\C$ is the free category with sums generated by $\C$.
Similarly, $\Famf\C$ is the free category with finite sums generated by $\C$.
Clearly, a universal functor $\C\to\Fam\C$ takes $A\in\C$ to the corresponding family $1\to \obj\C$, 
where $1$ is a terminal set.
A sum in $\Famf\C$ of $A:I\to\obj\C$ and $B:J\to\obj\C$ is given by the universally induced family $I+J \to \obj\C$,
while $0\to\obj\C$ is initial.
Thus, the fibration $\Fam\C \to \Set$ preserves sums.

Here are some instances of the family construction 
(some of which are given up to equivalence).
\begin{enumerate}
\item
If $\C = 1$, $\Fam\C = \Set$, and as a fibration it is the identity.
\item
If $\C = S$ is discrete, $\Fam\C$ is the domain of the discrete family fibration,
whose objects over $I$ are the families $s_i \, (i \in I)$ of elements of $S$ 
and arrows are given by reindexing.
Equivalently, $\Fam\C = \Set/S$.
\item
If $\C = \Set$, $\Fam\C$ is the arrow category $\Set^\to$. 
As a fibration, it is the codomain fibration $\Set^\to \to \Set$.
\item
More generally if $\C = \Set^\D$ then $\Fam\C= \Set^{\D'}$, where 
$\D'$ is obtained by $\D$ by adjoining a terminal object.
\item
If $\B = \core(\C)$ is the groupoid of isomorphisms of $\C$, 
$\Fam\B$ is the (domain of) the subfibration of $\Fam\C$ formed by its cartesian arrows.
\item
In particular, $\core(\Set)$ is the category of sets and bijections and
$\Fam(\core(\Set))$ is the wide subcategory $\Pb$ of $\Set^\to$, whose arrows are pullback squares.
\end{enumerate}

Recall that an object $A$ of a category $\D$ with sums is {\em connected}
if the representable functor $\D(A,-):\D\to\Set$ preserves sums.
\begin{prop} \label{fam1}
A category $\D$ is equivalent to $\Fam\C$ (respectively, $\Famf\C$) for some $\C$ 
if and only if any object of $\D$ is a (finite) sum of connected objects in $\D$.
\end{prop}
\pf
In one direction, any object of $\Fam\C$ is a sum of families indexed by a terminal set,
which are clearly the connected objects of $\Fam\C$. 
In the other direction, it is straightforward to see that, in that hypothesis, $\D \iso \Fam\C$, 
where $\C$ is the full subcategory of connected objects.
\epf

Let us denote by $\dCat$, $\sCat$ and $\sdCat$ the wide subcategories of $\Cat$
formed by discrete fibrations, sum preserving functors and sum preserving 
discrete fibrations, respectively.
\begin{prop} \label{fam2}
The adjunction between $\Cat$ and $\sCat$ given by the family construction restricts 
to an adjunction between $\dCat$ and $\sdCat$.
\end{prop}
\pf
Since, for a functor $F:\C \to \D$, the functor $\Fam F: \Fam\C \to \Fam\D$ takes
an arrow (over $f:I\to J$) $\al_i : A_i \to B_{fi} \, (i\in I)$ in $\Fam\C$ to the arrow
$F\al_i : FA_i \to FB_{fi} \, (i\in I)$ in $\Fam\D$, it is immediate to note that, 
if $F$ is a discrete fibration, then also $\Fam F$ is a discrete fibration.
\epf
\begin{prop} \label{fam3}
Suppose that $\E$ is a category with finite sums, that $F:\E\to\C$ is a sum preserving discrete fibration, 
and that $X\in\E$ is over a sum of a family $A_i\in\C$. 
Then $X$ is itself a sum of a family $X_i\in\E$ over $A_i$.
\end{prop}
\pf
Say, for instance, that $FX = A$ is a sum $A = B+C$ with injections $i:B\to FX$ and $j:C\to FX$.
Then, since $F$ is a discrete fibration, we have (unique) liftings $i':Y\to X$ and $j':Z\to X$.
Consider the sum $X' = Y+Z$ in $\E$, with injections $k:Y\to X'$ and $l:Z\to X'$.
By the hypothesis, $Fk:FY\to FX'$ and $Fl:FZ\to FX'$ exhibit $FX'$ as a sum: $FX' = FY+FZ = B+C$.
Since also $FX = B+C$, we have a comparison isomorphism $FX \iso FX'$ which lifts to an isomorphism $f:X\iso X'$
(since $F$ is a discrete fibration), with $fi' = k$ and $fj' = l$; thus also $i':Y\to X$ and $j':Z\to X$ exhibit $X$ as a sum. 
\epf
\begin{prop} \label{fam4}
If $F:\E\to\Fam\C$ is a sum preserving discrete fibration, then $\E = \Fam\D$ 
and $F=\Fam G$ for a discrete fibration $G:\D\to\C$
\end{prop}
\pf
This follows directly from Propositions \ref{fam1} and \ref{fam3}.
\epf


\subsection{Double categories and double props}
\label{Dc}\qq 
We assume the basic notions about double categories and functors.
In the present section we just fix notations and recall some instances and properties 
that will be useful in the sequel.

A double category $\dA$ is a category in $\Cat$, with its underlying
graph in $\Cat$:
\eq 
\label{dc2}
s : \dA_1 \to \dA_0 \qv t : \dA_1 \to \dA_0 
\eeq
The functor $\obj:\Cat\to\Set$ takes a double category $\dA$
to an ordinary category $\dA_l$.

We adopt the convention of calling \"tight arrows" the arrows of $\dA_0$,
while the \"loose arrows" are the objects of $\dA_1$, which are also the arrows of $\dA_l$.
The \"cells" or \"squares" are the arrows of $\dA_1$.
Despite the \"asymmetric" terminology adopted, the double categories considered here are strict: 
the \"external" composition $\dA_1\times_{\dA_0}\dA_1\to \dA_1$ 
(which makes the graph (\ref{dc2}) a category in $\Cat$) is strictly associative.
By the well-known intrinsic symmetry of strict double categories, the cells are also the arrows 
of another category whose objects are tight arrows.

We say that a double category $\dA$ has (finite) products or sums 
if both $\dA_0$ and $\dA_1$ have (finite) products or sums
and the source and target functors (\ref{dc2}) preserve them.
(These double categories are called \"precartesian" in \cite{patterson}.)
Clearly, $\dA$ has sums if and only if its opposite  $\dA\op$ 
(obtained by taking the opposite of $\dA_0$ and of $\dA_1$) has products.

By a (double) functor $F:\dA \to \dB$, with components $F_0:\dA_0 \to \dB_0$ 
and $F_1:\dA_1 \to \dB_1$ we intend a strict one.
It also has a \"loose part" $F_l:\dA_l \to \dB_l$.
We say that $F:\dA \to \dB$ preserves products or sums if both its components do so.

For a category $\C$, we denote by $\dSq\C$ the double category of commutative squares in $\C$, 
and by $\dPb\C$ the double subcategory of pullback squares.
If $\C = \Set_f$ we denote the latter simply by $\dPb$.
We denote by $\dBij$ the double full subcategory of $\dPb$ 
(or, equivalently, of $\dSq\Set_f$) whose loose arrows are bijections. 

Other double categories which we will use in the sequel are obtained 
by adding some kind of structure to $\dPb$.
For instance, we denote by $\dTot$ the double category obtained by endowing 
each fiber of the loose arrows of $\dPb$ with a total order.
The external composition is the obvious one, obtained by \"glueing" total orders along a total order.
The double category of sections $\dSec$ is obtained by endowing the loose arrows of $\dPb$ with a section.

Note that all these double categories have $\dA_0 = \Set_f$ and have finite sums,
that is they are \"double props" in the following sense:
\begin{definition}
\label{defprop}
A {\em double prop} is a double category $\dP$ with finite sums whose
tight category $\dP_0$ is the category $\Set_f$ of finite sets.
\end{definition}
The terminology is motivated by the fact that such a double category indeed gives a
prop in the classical sense, namely (a skeleton of) its loose part $\dP_l$.
Indeed, $\dP_l$ is a monoidal category with the tensor product inherited by the sum 
in $\dP_0 = \Set_f$ (for the objects) and in $\dP_1$ (for the arrows).
On the other hand, the term \"double prop" is slightly misleading,
since it may suggest that it is to a prop what a double category is to a category,
which clearly is not the case.
Rather, the idea is that the notion of double prop may be considered the \"right" notion of prop.
Indeed, the definition is rather natural and most of the main instances of classical props are 
in fact traces of double props.

We will also consider the many sorted version of double prop:
\begin{definition}
\label{defprop2}
A {\em colored double prop} is a double category $\dP$ with finite sums whose
tight category $\dP_0$ is the category $\Fam S$ for a set $S$ of \"colors" or \"objects".
\end{definition}


\subsection{Species as graphs in $\Cat$}
\label{Spe} \qq
The notion of species was introduced in \cite{joyal} as a categorical version of the concept of structure,
allowing to develop in an invariant and effective way several combinatorial notions and properties.

It is known that there are connections between species and operads. 
In particular, it is not new the idea that an operad is nothing but a monoid in a suitable monoidal 
category of species (see \cite{kelly}; see also \cite{fiore} for the colored case).
Here we present our version of this fact and, in doing so, we also hopefully enlighten some 
aspects of the theory of species.

There are at least two reasons to present these connections.
On the one hand, we will see that the main instances of operads or multicategories 
can in fact be usefully seen as basic species with a natural monoid structure. 
On the other hand, species are most naturally defined in an intrinsically unbiased way: 
the base category is that of sets and bijections, and not one of its skeletons.
Hence, taking wholeheartedly the view of multicategories as \"multiplicative species",
we are led toward the unbiased notion of the former which is advocated here.  

Recall that the category of (finitary) species can be defined as the category 
\eq  \label{spe1}
\Sp = \dCat/\core(\Set_f)
\eeq 
of discrete (op)fibrations over $\core(\Set_f)$, the category of finite sets and bijections.
Thus a species is a discrete fibration $F:\D\to\core(\Set_f)$, where the objects of $\D$ over a set $I$
are to be thought of as the structures on $I$ of that species, and the arrows of $\D$ over a mapping $f:I'\to I$
say how the structures over $I$ are \"transported" along $f$ to structures over $I'$.

Since $\Pb = \Fam(\core(\Set_f))$ (see Section~\ref{Fam}), 
by Proposition \ref{fam4} to give a discrete fibration $F:\D\to\core(\Set_f)$ is equivalent to give  
a category $\E$ with finite sums along with a sum preserving discrete fibration $G:\E\to\Pb$
(with $\E = \Fam\D$ and $G=\Fam F$).
That is, the category of species can be also defined as the slice
\eq \label{spe2}
\Sp = \sdCat/\Pb
\eeq 
Now, $\Pb$ has a natural structure of double graph (meaning, a graph in $\Cat$),
inherited by its inclusion in $\Set_f^\to$
\eq \label{spe3}
s,t : \Pb \to \Set
\eeq 
so that any species $G:\E\to\Pb$ has a graph structure over it
\eq \label{spe4}
sG,tG : \E \to \Set
\eeq 
It is then natural to define colored species as follows.
Let $\DG$ be the category which has
\begin{itemize}
\item
as objects, double graphs $\dG$ with sums,
that is graphs $s,t : \dG_1 \to \dG_0$ in $\Cat$ such that $\dG_1$ and $\dG_0$ 
have finite sums preserved by $s$ and $t$; 
\item
as arrows, morphisms $G:\dG \to \dH$ of double graphs such that both the components
$G_1:\dG_1 \to \dH_1$ and $G_0:\dG_0 \to \dH_0$ are sum preserving discrete fibrations.
\end{itemize}
Then the category of colored species is the slice
\eq \label{spe5}
\colSp = \DG/\dPb
\eeq 
where $\dPb$ is of course the underlying graph of the corresponding multicategory.

Note that if $G:\dG \to \dPb$ is a colored species, 
\eq
\label{spe7}
\xymatrix@R=1.5pc@C=1.5pc{
\dG_1 \ar[dd]_{G_1}\ar[rr]<3pt>^s\ar[rr]<-2pt>_t & & \dG_0 \ar[dd]^{G_0} \\
& & \\
\Pb  \ar[rr]<3pt>^s\ar[rr]<-2pt>_t                 & & \Set_f       }
\eeq
then the component $G_0:\dG_0\to\Set_f$, 
being a sum preserving discrete fibration over $\Set_f = \Famf(1)$
is (again by Proposition \ref{fam4}) simply the family fibration on a set $S$,
which we can think of as the set of \"objects" or \"sorts" of the species.
Intuitively, each structure over $I$ has a codomain in $S$, and a domain which is a family $I\to S$.
Colored species can be considered as the unbiased version of \"colored symmetric sequences".

We can then identify the category $\Sp$ as the full subcategory of $\colSp$ spanned by \"untyped" 
or \"single-sorted" species, that is those for which $G_0$ is the identity (or an isomorphism).


\subsection{Changing the base}
\label{Chba} \qq
Now it is clear that the \"relative" approach to multicategories, sketched in Section \ref{Bas}, 
can be considered also at the level of species 
(that is with double graphs in place of double categories).
So, if $\dG \in \DG$, the category of $\dG$-species is the slice
\eq  \label{spe6}
\dG\-\Sp = \DG/\dG
\eeq 
As above, we say that $G:\dG'\to\dG$ is untyped over $\dG$ if 
$G_0:\dG_0'\to\dG_0$ of $G$ is the identity (or an isomorphism).

Let us consider some examples of $\dG$-species, in which $\dG$ is itself 
(the domain of) a species $\dG\to\dPb$.
\begin{enumerate}
\item
Consider the species \"characteristic of singletons", which has 
just one structure over the terminal sets and none over the other sets.
From our point of view, it corresponds to the following morphism $\dBij \to \dPb$ in $\DG$:
\eq
\label{spe8}
\xymatrix@R=1.5pc@C=1.5pc{
\Bij \ar[dd]\ar[rr]<3pt>^s\ar[rr]<-2pt>_t & & \Set_f \ar[dd]^{\id} \\
& & \\
\Pb  \ar[rr]<3pt>^s\ar[rr]<-2pt>_t                 & & \Set_f       }
\eeq
where $\Bij = \dBij_1$ (see Section \ref{Dc}) is the full subcategory
of $\Pb = \dPb_1$ or of $\Set_f^\to$ whose objects are bijective mappings
(that is, mappings with a \"characteristic of singletons" structure on each fiber).
Then the \"untyped" $\dG$-species are sets, while the typed ones are graphs.
\item
Consider the species of total orders.
From our point of view, it corresponds to the following morphism $\dTot \to \dPb$ in $\DG$:
\eq
\label{spe9}
\xymatrix@R=1.5pc@C=1.5pc{
\Tot \ar[dd]\ar[rr]<3pt>^s\ar[rr]<-2pt>_t & & \Set_f \ar[dd]^{\id} \\
& & \\
\Pb  \ar[rr]<3pt>^s\ar[rr]<-2pt>_t                 & & \Set_f       }
\eeq
where $\Tot = \dTot_1$ is the category over $\Pb$ whose objects are mappings 
with totally ordered fibers.
Then the untyped $\dTot$-species are \"totally ordered" species,
that is those species whose structures have a total order on their underlying set. 
The $S$-typed ones may be called \"planar multigraphs": $S$ is the set of objects and each arrow has a codomain in $S$, 
while its domain is a totally ordered family in $S$ (that is, a family $I\to S$ with a total order on $I$).
\item
Consider the species elements, whose structures over a set $I$ are the elements of $I$.
Note that, as a discrete opfibration on $\core(\Set_f)$, it corresponds to the inclusion functor $\core(\Set_f) \to \Set$.
From our point of view, it corresponds to the following morphism $\dSec \to \dPb$ in $\DG$:
\eq
\label{spe8b}
\xymatrix@R=1.5pc@C=1.5pc{
\Sec \ar[dd]\ar[rr]<3pt>^s\ar[rr]<-2pt>_t & & \Set_f \ar[dd]^{\id} \\
& & \\
\Pb  \ar[rr]<3pt>^s\ar[rr]<-2pt>_t                 & & \Set_f       }
\eeq
where $\Sec = \dSec_1$ (see Section \ref{Dc}) is the category whose objects are mappings with a section,
that is pairs $(f,s)$ such that $fs$ is the identity and whose arrows are pullbacks (with the pulled back sections).
Then the \"untyped" $\dSec$-species are \"pointed" species,
that is those species whose structures have a selected element of their underlying set. 
\item
For any set of \"labels" $L$, one has the species $\dG_L$ whose structures over $I$ are 
the $L$-labelings of $I$, that is the mappings $I\to L$ in $L^I$ 
(this is again the discrete family fibration on $L$, restricted to bijections).
Then the untyped $\dG_L$-species are species labeled by $L$.

If $L$ is the set of arrows of a graph $\L$ then we can define the following $\dG_L$-species $\dG_\L$,
typed by the objects of $\L$: 
the structures over $I$ are families $I\to L$ of concurrent arrows in $\L$, with domain and codomain given by
the graph structure of $\L$.
This example will be important in the sequel, when from a category $\L$ 
we so get a symmetric multicategory (see Section \ref{Ex}).

\end{enumerate}


\subsection{Adding structure cofreely}
\label{Adstr} \qq
In general, if $\dG\to\dPb$ is a species and $G:\dG'\to\dG$ is a $\dG$-species, composing with $G$ we
get a \"forgetful" functor $G_!: \DG/\dG' \to\DG/\dG$ from $\dG'$-species to $\dG$-species.
Since $\dCat$ has pullbacks, $G_!$ has a right adjoint given by pulling back along $G$:
\eq
G_!\adj G^*:\dG\-\Sp \to\dG'\-\Sp
\eeq

For instance, considering the species $G:\dTot\to \dPb$ of total orders, we get an adjunction between species and ordered species:
\eq
G_!\adj G^*:\Sp \to\dTot\-\Sp
\eeq
The left adjoint simply \"forgets" the order, while the right adjoint $G^*$ (being the pullback) multiplies the fibers:
for any species $H:\dH\to \dPb$, the structures of the ordered species $G^*H$ are obtained by pairing 
each $H$-structure over $I$ with all possible orders on $I$.
Similarly, we have
\eq
G_!\adj G^*:\Sp \to\dBij\-\Sp
\eeq
where the right adjoint just retains the structures over the terminal sets.


\section{Unbiased symmetric muticategories}
\label{USM}

In this section, we return to our definition of unbiased symmetric muticategory:
a double category $\dM$ with finite sums with a sum preserving discrete fibration 
\( M:\dM\to\dPb \) to the double category of pullback squares in finite sets.


\subsection{Unbiased symmetric multicategories}
\label{Usm} \qq
Recall that in Section \ref{Bas} we considered the category $\DP$ 
whose objects are colored double props (see Section \ref{Dc}) and whose arrows are sum preserving discrete fibrations.
Then we defined the category $\sMlt$ of (unbiased) symmetric multicategories 
as the slice of $\DP$ over the double prop $\dPb$:
\eq \label{usm1}
\sMlt = \DP / \dPb 
\eeq


\subsection{Symmetric multicategories as monoids on species}
\label{Smon} \qq
The point of view of species as double graphs of Section \ref{Spe} 
has the advantage of rendering natural the consideration of a category structure on them,
allowing us to consider symmetric multicategories as monoid structures on species.
The idea is that the underlying colored species gives the arrows with their symmetry, 
while the monoid structure gives the composition of arrows. 

In fact, the wide subcategory $\Pb$ of $\Set_f^\to$ inherits from it the structure 
of a double category $\dPb$, a category in $\Cat$. 
Thus, $\dPb$ is a monoid in the bicategory $\Span\Cat$ or, equivalently,
on the monoidal category of double graphs which have $\Set_f$ has category of objects.
Then the category of graphs over 
\eq \label{smon1}
s,t : \Pb \to \Set_f
\eeq 
inherit a monoidal structure, as a consequence of the following general fact:
\begin{prop} 
\label{smon3}
If $\C$ is a monoidal category and $(X,\mu,\eta)$ is a monoid in $\C$,
then $\C/X$ has a natural monoidal structure such that the monoids in $\C/X$ are 
the monoids in $\C$ over $X$: $\Mon(\C/X) = (\Mon\,\C) / X$.
\end{prop}
\pf
The tensor of $f:A\to X$ and $g:B\to X$ is given by 
$f\otimes' g = \mu\circ (f\otimes g):A\otimes B\to X$.
The rest is easily checked.
\epf

Recall that in Section \ref{Spe} we have defined a colored species as 
a double graph $\dG$ with sums, along with a sum preserving discrete fibration
$G:\dG \to \dPb$.
It is then a commuting diagram in $\sCat$, where $G_0$ and $G_1$ are discrete fibrations
(so that $\dG_0 = \Fam S$, for a set $S$ of \"objects"):
\eq
\label{smon2}
\xymatrix@R=1.5pc@C=1.5pc{
\dG_1 \ar[dd]_{G_1}\ar[rr]<3pt>^s\ar[rr]<-2pt>_t & & \dG_0 \ar[dd]^{G_0} \\
& & \\
\Pb  \ar[rr]<3pt>^s\ar[rr]<-2pt>_t                 & & \Set_f      }
\eeq
Then, 
\begin{prop} 
\label{smon4}
An unbiased symmetric multicategory amounts to a monoid structure on a colored species, 
for the monoidal structure on $\colSp$ inherited by the double category structure on $\dPb$.
\end{prop}
\pf
Indeed, by definition, an unbiased symmetric multicategory is a double category structure 
over a colored species $\dG$ as in (\ref{smon2}). 
Then the result follows by Proposition \ref{smon3}.
\epf
If $G_0$ is an isomorphism (that is the species $\dG$ is \"untyped", or, equivalently, typed by a singleton),
the corresponding multicategory has a single object and we get \"unbiased operads".


\subsection{Examples}
\label{Ex} \qq
We now view (or review) some important instances of symmetric multicategories,
which are most naturally seen as species with a monoid structure.

Of course
\eq
\label{ex1}
\xymatrix@R=1.5pc@C=1.5pc{
\Pb \ar[dd]_{\id}\ar[rr]<3pt>^s\ar[rr]<-2pt>_t & & \Set_f \ar[dd]^{\id} \\
& & \\
\Pb  \ar[rr]<3pt>^s\ar[rr]<-2pt>_t                 & & \Set_f      }
\eeq
is the terminal species, the species with just one structure on any set.
It is also the terminal symmetric multicategory $1\t$ (we will soon explain that notation).

The species characteristic of singletons (see Section \ref{Spe}) 
\eq
\xymatrix@R=1.5pc@C=1.5pc{
\Bij \ar[dd]\ar[rr]<3pt>^s\ar[rr]<-2pt>_t & & \Set_f \ar[dd]^{\id} \\
& & \\
\Pb  \ar[rr]<3pt>^s\ar[rr]<-2pt>_t                 & & \Set_f       }
\eeq
where $\Bij = \dBij_1$, is a symmetric multicategory
given by the double category structure of $\dBij$.

The species of total orders (see Section \ref{Spe}) 
\eq
\xymatrix@R=1.5pc@C=1.5pc{
\Tot \ar[dd]\ar[rr]<3pt>^s\ar[rr]<-2pt>_t & & \Set_f \ar[dd]^{\id} \\
& & \\
\Pb  \ar[rr]<3pt>^s\ar[rr]<-2pt>_t                 & & \Set_f       }
\eeq
where $\Tot = \dTot_1$, is a symmetric multicategory given by the double category structure of $\dTot$ (see Section \ref{Dc}).

The species elements (see Section \ref{Spe}) 
\eq
\label{ex4}
\xymatrix@R=1.5pc@C=1.5pc{
\Sec \ar[dd]\ar[rr]<3pt>^s\ar[rr]<-2pt>_t & & \Set_f \ar[dd]^{\id} \\
& & \\
\Pb  \ar[rr]<3pt>^s\ar[rr]<-2pt>_t                 & & \Set_f       }
\eeq
where $\Sec = \dSec_1$, is a symmetric multicategory
given by the double category structure of $\dSec$.

If $\C$ is a category, we get a symmetric multicategory $M_\C:\dC\t \to\dPb$
\eq
\label{seq}
\xymatrix@R=1.5pc@C=1.5pc{
(\dC\t)_1 \ar[dd]\ar[rr]<3pt>^s\ar[rr]<-2pt>_t & & (\dC\t)_0 \ar[dd]^{\id} \\
& & \\
\Pb  \ar[rr]<3pt>^s\ar[rr]<-2pt>_t                 & & \Set_f       }
\eeq
where $(\dC\t)_0$ is the discrete family fibration $\Fam(\obj\C)$ 
on the objects of $\C$, while $(\dC\t)_1$ has 
\begin{itemize}
\item
as objects, the arrows of $\Fam\C$ (which are then the loose arrows of $\dC\t$);
\item
as arrows, the commuting squares in $\Fam\C$ over pullback squares in $\Set_f$,
whose horizontal (say) arrows are cartesian lifting of those of the underlying 
pullback.
\end{itemize}
Note that the underlying double graph of $\dC\t$ is the graph $\dG_\C$ 
labeled by the arrows of $\C$ (see the last example in section \ref{Spe}).

\subsection{Cocartesian multicategories}
\label{Coca} \qq
The symmetric multicategories $\dC\t$ are presented in \cite{hermida} under the 
(slightly unwieldy) name \"discrete cocone" multicategories, since their arrows are sequences 
(or families, in the unbiased approach) of concurrent arrows in $\C$; 
the notation $\dC\t$ is also borrowed from that paper.

They are studied in detail in \cite{pisani} under the name of \"sequential multicategories", 
a rather unhappy choice especially for the present unbiased context. 

A far better name, suggested by Nathanael Arkor, is \"cocartesian multicategories",
since they are representable (that is, monoidal) exactly when $\C$ has finite sums
(see Section \ref{Smc}).
This would nicely parallel the case of cartesian multicategories (see Section \ref{CART}),
which are representable exactly when the underlying category has finite products.


\subsection{Changing the base}
\label{Bas2}\qq
As sketched in the introduction, we can usefully relativize the notion
of unbiased multicategory by using any double prop $\dP$ in place of $\dPb$ as the indexing base,
so getting the category of $\dP$-multicategories:  
\eq 
\dP\-\Mlt = \DP / \dP 
\eeq 
so that
\eq 
\dPb\-\Mlt = \sMlt 
\eeq 
For instance, 
\eq 
\dBij\-\Mlt = \Cat 
\eeq 
and 
\eq 
\dTot\-\Mlt = \Mlt 
\eeq 


As explained in the case of species in Section \ref{Adstr},
if $\dP\to\dPb$ is a symmetric multicategory and $P:\dP'\to\dP$ is a $\dP$-multicategory, composing with $P$ we
get a \"forgetful" functor $P_!: \DP/\dP' \to\DP/\dP$ from $\dP'$-multicategories to $\dP$-multicategories,
which has a right adjoint given by pullbacks:
\eq
P_!\adj P^*:\dP\-\Mlt \to\dP'\-\Mlt
\eeq

For instance, considering the \"associative operad" $T:\dTot\to \dPb$, we get 
an adjunction between symmetric and plain multicategories:
\eq
T_!\adj T^*:\sMlt \to\Mlt
\eeq
For a plain multicategory $M'$, the arrows of $T_!M'$ are obtained by \"forgetting" the order on the domain of the arrows of $M'$, 
while, for a symmetric multicategory $M$, the arrows of $T^*M$ are obtained by attaching any possible total order on the domain
of each arrow of $M$. Similarly, considering the operad $U:\dBij \to \dPb$ (with exactly one object and one arrow), we get 
an adjunction between symmetric multicategories and categories:
\eq
U_!\adj U^*:\sMlt \to\Cat
\eeq
where the left adjoint takes a category to the corresponding unary multicategory,
while the right adjoint takes a multicategory to the category formed by its unary arrows.

\begin{remark}
It is not difficult to see that a symmetric multicategory $P:\dP \to \dPb$ admits a plain multicategory structure
$P':\dP\to\dTot$ (with $P=T\circ P'$) if and only if, as a species, it is \"flat" in the sense of \cite{bergeron} 
(see also \cite{gambino}) that is the associated group action is free (or, equivalently, each structure is \"rigid").
So, while classically one considers symmetric multicategories as plain multicategories
with further structure, namely the action on arrows associated to the permutations of the domain,
we here rather consider plain multicategories as symmetric multicategories with a property, namely the rigidity of arrows.

Likewise, one may ask whether an arrow in a symmetric multicategory $M$, 
say $\al: A,A \to B$ (with domain $2\to 1\to \obj M$) is \"symmetric" (or \"commutative")  that is, 
whether its reindexing along the non-identity mapping $2\to 2$ is $\al$ again; 
however, classically, the same question cannot be posed for an arrow in a plain multicategory, 
since we have not reindexing of arrows therein.
From our point of view, we would rather say that an arrow in a plain multicategory cannot have symmetries because 
the reindexing is free (arrows are rigid).
\end{remark}


\subsection{Unbiased monoidal categories}
\label{Smc}\qq
As already sketched in the Introduction, we can consider the unbiased versions of monoidal category 
and of monoids, relative to any base double prop $\dP$.
\begin{definition} 
An {\em unbiased $\dP$-monoidal category} is a $\dP$-multicategory $M: \dM\to\dP$ 
such that $M_l: \dM_l\to \dP_l$, the loose part of $M$, is an opfibration whose opcartesian arrows 
are stable with respect to reindexing.
We also say that the $\dP$-multicategory $M$ is {\em representable}.
\end{definition} 
In fact, when $\dP = \dPb$, this is equivalent to the usual condition for expressing the existence of tensor products 
in a multicategory via the universal form of representability (see \cite{hermida}, \cite{leinster} and \cite{pisani}).
\begin{definition} 
An {\em unbiased $\dP$-monoid} (or {\em $\dP$-algebra}) is a $\dP$-multicategory $M: \dM\to\dP$ 
such that $M_l: \dM_l\to \dP_l$, the loose part of $M$, is a discrete opfibration.
\end{definition} 
This (when $\dP = \dPb$) is an expression of the fact that monoids arise as the discrete multicategories
(or the discrete monoidal categories).

Then, for any double prop $\dP$, we have the category of $\dP$-monoidal categories
\eq
\dP\-\MonCat \subseteq \dP\-\Mlt
\eeq 
and  the category of $\dP$-monoids
\eq
\dP\-\Mon \subseteq \dP\-\MonCat
\eeq

Thus, for the standard indexing $\dP = \dPb$ we have
\eq 
\dPb\-\MonCat = \sMonCat 
\eeq 
the category of (unbiased) symmetric monoidal categories, and
\eq 
\dPb\-\Mon = \cMon
\eeq 
the category of (unbiased) commutative monoids.
More generally, for $\dP = \dC\t$ we get
\eq
\dC\t\-\MonCat = \sMonCat^\C
\eeq 
that is, a $\dC\t$-monoidal category corresponds to a functor from $\C$ to 
the category of monoidal categories and strong monoidal functors.
Similarly,
\eq
\dC\t\-\Mon = \cMon^\C
\eeq
For the indexing $\dP = \dTot$ we have
\eq 
\dTot\-\MonCat = \MonCat 
\eeq 
the category of (unbiased) monoidal categories, and 
\eq 
\dTot\-\Mon = \Mon
\eeq 
the category of (unbiased) monoids.
For $\dP = \dBij$ we get
\eq
\dBij\-\MonCat = \Cat
\eeq 
and 
\eq
\label{bijmon}
\dBij\-\Mon = \Set
\eeq

If $M:\dM \to \dPb$ is a symmetric multicategory, then $\dM\-\Mon$
is the usual category of $M$-algebras,
that is, of functors to the symmetric multicategory of sets and several variable mappings.
The same holds for a plain multicategory $M:\dM \to \dTot$.

If $\C$ is a category with finite sums, the symmetric multicategory $M_\C:\dC\t \to\dPb$  
(see (\ref{seq}) in Section \ref{Ex}) is monoidal.
Indeed, its loose part is the family fibration on $\C$, $\Fam\C \to \Set_f$, and the latter is a bifibration
exactly when $\C$ has finite sums.

Other instances of unbiased monoidal categories are given in Section \ref{Unmon}.


\subsection{The Burnside monoid of a monoidal category}
\label{Bmon} \qq

By taking the isomorphism classes of a (symmetric) monoidal category, 
one gets a (commutative) monoid. 

In the present context, this can be generalized as follows:
\begin{prop}
By taking the isomorphism classes of the loose part $\dM_l$ of a $\dP$-monoidal category $M:\dM\to\dP$, 
one gets a $\dP$-monoid.
\end{prop}
\pf
Indeed, this follows easily from the stability condition on opcartesian arrows.
\epf


\subsection{Commutative monoids}
\label{Cm}\qq 
We have just defined an (unbiased) commutative monoid $M$ as a \"discrete" symmetric multicategory
$M:\dM\to\dPb$, that is one whose loose part $M_l$ is a discrete opfibration.
Nevertheless, we present here directly this special case, since it has its own importance and
may provide some further motivation for the more general related notions.

It is natural to consider a commutative and associative operation on a set $M$
as providing a way to \"add" any $I$-indexed family of elements of $M$ along
any mapping $f:I\to J$ to give a $J$-indexed family.
This can be taken as the starting point of an unbiased notion of commutative monoid.

Indeed, given a commutative monoid $M$, we can define a functor $\Set_f \to \Set$
which is the family functor $I \mapsto M^I$ on objects and acts on arrows by adding up families:
if $x = x_i \, (i\in I)$ is a family in $M^I$, we define its \"sum" $f_! x$ along a mapping $f:I\to J$ as 
the family $(f_!x)_j = \Sigma_{fi=j} x_i$ in $M^J$.

Furthermore, this functor has the following compatibility condition with the reindexing functor
$\Set_f\op \to \Set$ (that is the functor represented by $M$, which is $I \mapsto M^I$ on objects and 
is given by composition on arrows: $(f^*x)_i = x_{fi}$).
For any pullback in $\Set_f$ the corresponding right hand square below commutes:
\eq
\label{go1}
\xymatrix@R=1.3pc@C=1.3pc{
I \ar[dd]_f\ar[rr]^k  & & L \ar[dd]^g \\
& \pb & \\
J  \ar[rr]_l                 & & K       }
\qq\qq
\xymatrix@R=1.8pc@C=1.8pc{
M^I \ar[dd]_{f_!}  & & M^L \ar[ll]_{k^*}\ar[dd]^{g_!} \\
&  & \\
M^J                  & & \ar[ll]^{l^*} M^K       }
\eeq
This compatibility condition expresses, when $l:J\to K$ is a bijection, the fact that the result of the addition 
does not depend on the particular indexing. On the other hand, when $l:1\to K$ picks an element, 
it says also that the result of adding an $L$-indexed family of elements of $M$ along a mapping $g:L\to K$ 
can be obtained by adding separately the families indexed by the fibers of $g$;
that is, addition along general mappings are determined by addition along mappings $L\to 1$.
 
Passing to a fibrational perspective, given a monoid $M$, we thus have a discrete fibration
(reindexing) and a discrete opfibration (addition) over $\Set_f$ which have the same objects 
(the families of elements of $M$). 
Furthermore, the discrete fibration and the discrete opfibration are compatible 
in the sense that transporting an $L$-family $x = x_i \, (i\in I)$ along the two possible 
paths over a pullback, we end up with the same $J$-family:
\eq
\label{go2}
\xymatrix@R=1.3pc@C=1.3pc{
I \ar[dd]_f\ar[rr]^k  & & L \ar[dd]^g \\
& \pb & \\
J  \ar[rr]_l                 & & K       }
\qq\qq
\xymatrix@R=1.8pc@C=1.8pc{
k^*x \ar@{|->}[dd]  & & x \ar@{|->}[ll]\ar@{|->}[dd] \\
&  & \\
f_!k^*x = l^*g_!x             & & \ar@{|->}[ll] g_!x      }
\eeq

We can then define a double category $\dM$ as follows:
\begin{itemize}
\item
the objects of $\dM$ are finite families $x = x_i \, (i\in I)$ of elements of $M$;
\item
the tight arrows of $\dM$ are those of the reindexing fibration;
\item
the loose arrows of $\dM$ are those of the addition opfibration;
\item
the cells of $\dM$ are the squares as in (\ref{go2}) over pullback squares in $\Set_f$.
\end{itemize}
and an obvious functor $M:\dM\to \dPb$ such that:
\begin{itemize}
\item
$M_0$ is a sum preserving discrete fibration and $M_l$ is a discrete opfibration;
\item
$M_1$ is a discrete fibration, because of (\ref{go2});
\item
$M_1$ is sum preserving. 
Indeed, the objects of $\dPb_1$, that is the mappings $f:I\to J$, are sums 
of their fibers $f^{-1}\{j\} \to 1$, for $j:1\to J$,
while the corresponding addition mapping $f_!:M^I \to M^J$ is the sum 
of the additions along the fibers $M^{f^{-1}\{j\}} \to M$.
\end{itemize}

Conversely, it easy to see that any double functor $M:\dM\to \dPb$ such that
$M_0$ and $M_1$ are sum preserving discrete fibrations 
and $M_l$ is a discrete opfibration, comes in such a way from a unique commutative monoid $M$.


\subsection{Infinitary multicategories}
\label{Inf} \qq
It is clear that in all we said the category $\Set_f$ can be replaced by the category $\Set$.
That is, we don't really need the operations to have finite families as domain.
Thus, an advantage of the present approach is that of providing a natural definition of
(unbiased) {\em infinitary} (symmetric) multicategory,
{\em infinitary} (symmetric) monoidal category and {\em infinitary} (commutative) monoid.

For instance, a category with sums or products is an infinitary symmetric monoidal category
and by taking its isomorphism classes we get an infinitary commutative monoid
(see Section \ref{Bmon}).

Another approach to infinitary symmetric monoidal categories is presented in \cite{jane}.
Even if the two resulting notions are likely equivalent, it is apparent that the unbiased point of view
yields a far more natural definition than the skeletal one.



\section{Unbiased cartesian muticategories}
\label{CART}

In this section, we define a monad $(-)^\cart$ on $\sMlt$, whose algebras are (unbiased)
cartesian multicategories.
That we so indeed capture the right notion of cartesian multicategory
is confirmed by the fact that one can prove therein the equivalence between 
tensor products (representability), universal products and algebraic products
(see Proposition \ref{maincart}).

Cartesian multicategories can be considered as a way to encode algebraic theories.
With respect to Lawvere theories, they are usually simpler since the products are in principle only \"virtual".
For example, if $R$ is a ring, the Lawvere theory for $R$-modules is given by {\em matrices} with entries in $R$,
while the associated cartesian multicategory is given by {\em sequences} 
(or, in the unbiased approach, {\em families}) with entries in $R$ (see example \ref{cartex2}).

\begin{remark}
Usually, cartesian multicategories are considered as a sort of variation of symmetric multicategories, 
where the action of bijective mappings on arrows is extended to an action of all mappings.
From our point of view, we have seen that the action of bijective mappings on the arrows 
of a symmetric multicategory is related to the {\em contravariant} reindexing of arrows 
along pullbacks.
We will show now that the action on arrows in a cartesian multicategory is instead 
given by a sort of {\em covariant} reindexing of arrows, related to the contravariant 
one by some kind of Beck-Chevalley and Frobenius laws
(see also Remark \ref{rmkfrob}.)
\end{remark}


It is well known that any full sub-multicategory of a cartesian monoidal category is an instance 
of cartesian multicategory.
Less known is the fact that also categories enriched in commutative monoids give instances of cartesian multicategories.
In fact, they can be {\em characterized} as the cocartesian (or \"sequential") cartesian multicategories 
(see \cite{pisani} and Section \ref{Cartco}).
In order to motivate the definition of cartesian multicategory, let us take a look to these main instances.
 
\begin{example}
\label{cartex1}
Given a mapping out of a product in $\Set$, one can \"duplicate or delete the variables" obtaining another mapping.
For instance, given $\al : A\tm B \tm A \to D$ one gets $\al': B\tm C \tm A \to D$ by posing
\[ 
\al'(b,c,a) = \al(a,b,a) \,\,.
\]
To be precise, if $Y$ is the family of sets $\{1,2,3\} \to\obj\Set$ given by $1\mapsto B,2\mapsto C,3\mapsto A$ 
and $f:\{1,2,3\} \to\{1,2,3\}$ is $1\mapsto 3,2\mapsto 1,3\mapsto 3$, then the domain $X=f^*Y$ of $\al$
is the family of sets $\{1,2,3\}\to\obj\Set$ given by $1\mapsto A,2\mapsto B,3\mapsto A$;
$\al' = f_!\al$ is obtained by precomposing $\al$ with the obvious map $\ov f : B\tm C \tm A \to A\tm B\tm A$.
(Note that as particular instances of $\ov f$ we find projections and diagonals.)
Of course, the same holds in any finite product category $\C$, or in any full sub-multicategory
of a cartesian monoidal multicategory.
\end{example}

\begin{example}
\label{cartex2}
If $\C$ is any additive category (or also any category enriched in commutative monoids), 
then there is a natural cartesian structure on the corresponding cocartesian (or \"sequential") multicategory 
$\dC\t$ (see \cite{pisani} and Section \ref{Coca}). 
In particular in the one object case, a ring $R$ gives a cartesian multicategory (operad) $\dR\t$ 
whose arrows are families of element of $R$.
If $\al:\{1,2,3\} \to R$ is the arrow $1\mapsto b,2\mapsto c,3\mapsto a$ in $\dM_R$, 
and $f:\{1,2,3\} \to\{1,2,3\}$ is $1\mapsto 3,2\mapsto 1,3\mapsto 3$ as in the example above,
then $f_!\al$ is obtained by summing up or by inserting $0$'s: $1\mapsto c,2\mapsto 0,3\mapsto b+a$.
\end{example}

\begin{remark}
\label{cartex3}
Thus, the notion of cartesian multicategory may be seen as a common generalization
of finite product categories and $\cMon$-enriched categories. 
In this perspective, one may wonder what are, for instance, the morphisms of cartesian multicategories
$\dR\t \to \Set$ between the instances just considered.
It turns out that these are exactly the $R$-modules (see \cite{pisani}). 
\end{remark}


\subsection{The classical approach to cartesian multicategories}
\label{Classic}

Though the explicit definitions of cartesian multicategory do not abound in the literature,
there seems to be a general consensus on the following requirements
(see \cite{pisani} and the references therein):
\begin{enumerate}
\item
There is a covariant reindexing $\al \mapsto f_!\al$  of maps along mappings.
\item
This reindexing is functorial.
\item
It is compatible with composition from above.
\item
It is compatible with composition from below.
\end{enumerate}

Let us exemplify graphically these conditions: 
\begin{enumerate}
\item
The covariant reindexing involves a reindexing of objects
along a mapping $f$ as in the above examples. 
We denote the corresponding tight arrow in $\dM$ as follows:
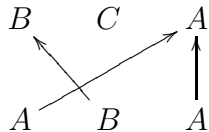
\[
\xymatrix@R=2pc@C=1.5pc{
B          & C  &   A                    \\
A \ar[urr] &        B  \ar[ul]          & A \ar[u] 
}
\]
Then, the covariant reindexing along $f$ associates an arrow $f_!\al:B,C,A\to D$ 
to any arrow $\al:A,B,A\to D$:
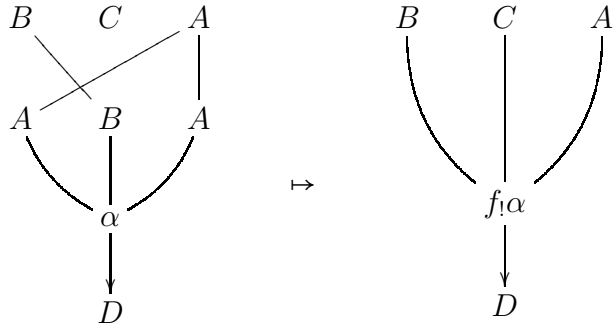
\[
\xymatrix@R=2pc@C=1.5pc{
B  \ar@{-}[dr]          & C  &   A  \ar@{-}[d]  \ar@{-}[dll]                   \\
A \ar@/_/@{-}[dr] &        B   \ar@{-}[d]          & A \ar@/^/@{-}[dl] \\
                         &  \al \ar[d]                           &                           \\
                         &  D                                       &
}
\xymatrix@R=2pc@C=1.5pc{
              &  &                       \\ \\
              & \mapsto & 
}
\xymatrix@R=2pc@C=1.5pc{
B  \ar@/_1pc/@{-}[ddr]     &   C \ar@{-}[dd]   &   A \ar@/^1pc/@{-}[ddl]                     \\
                                                 &                                &                           \\
                                                 &  f_!\al \ar[d]       &                            \\
                                                 &  D                            &
}
\]
Note that it can be seen as a way to evaluate \"spans" in $\dM$, with a loose leg and a tight leg
(see Section \ref{Cart}):
\eq 
\xymatrix@R=1.3pc@C=1pc{
& & f^*X \ar[ddll]_\al\ar@{-->}[ddrr]^f  & &  \\ 
& & & &  \\
Y & &                  & & X      } 
\xymatrix@R=2pc@C=1.5pc{
              &  &                       \\
              & \mapsto & 
}
\xymatrix@R=1.3pc@C=1pc{
  \\ 
& & & &  \\
Y & &                  & & X\ar@{..>}[llll]_{f_!\al}     } 
\eeq

\item
Functoriality of reindexing, \( f_!(g_!\al) = (f g)_!\al \),
may be illustrated by the fact that diagrams such as this one can be read unambiguously: 
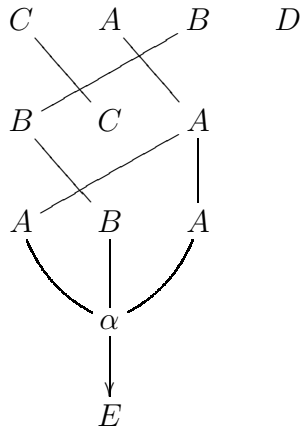
\[
\xymatrix@R=2pc@C=1.5pc{
C  \ar@{-}[dr]          & A \ar@{-}[dr]  &   B \ar@{-}[dll]   & D               \\
B  \ar@{-}[dr]          & C  &   A  \ar@{-}[d]  \ar@{-}[dll]                   \\
A \ar@/_/@{-}[dr] &        B   \ar@{-}[d]          & A \ar@/^/@{-}[dl] \\
                         &  \al \ar[d]                           &                           \\
                         &  E                                       &
}
\]

\item
Compatibility with composition from above
may be illustrated by the fact that diagrams such as this one can be read unambiguously: 
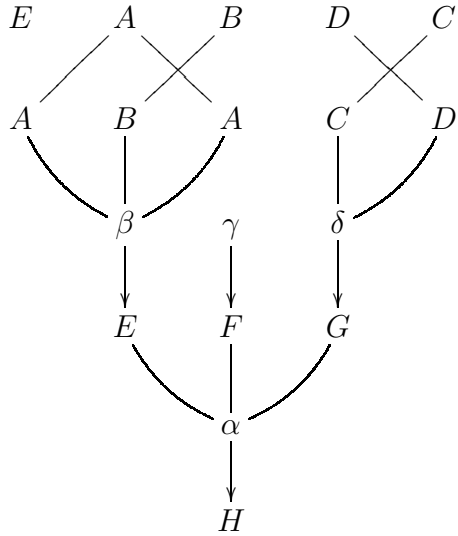
\[
\xymatrix@R=2pc@C=2pc{
E                & A \ar@{-}[dl]\ar@{-}[dr] &   B  \ar@{-}[dl]   & D\ar@{-}[dr]  & C\ar@{-}[dl]           \\
A  \ar@/_/@{-}[dr]                & B \ar@{-}[d] &   A  \ar@/^/@{-}[dl]   & C\ar@{-}[d]  & D\ar@/^/@{-}[dl]           \\
                 & \beta \ar[d]              &       \gamma \ar[d]                   & \delta \ar[d]                            \\
                 & E \ar@{-}@/_/[dr]    &          F   \ar@{-}[d]           & G \ar@{-}@/^/[dl]                            \\
                 &                               &  \al \ar[d]                              &      &                    \\
                 &        &  H                                          &      &
}
\]

\item
Compatibility with composition from below may be illustrated by the fact that, in a diagram such as this one, 
the covariant reindexing can be \"lifted" on the top 
(sometimes this is called, for obvious reasons, the \"combing" condition): 
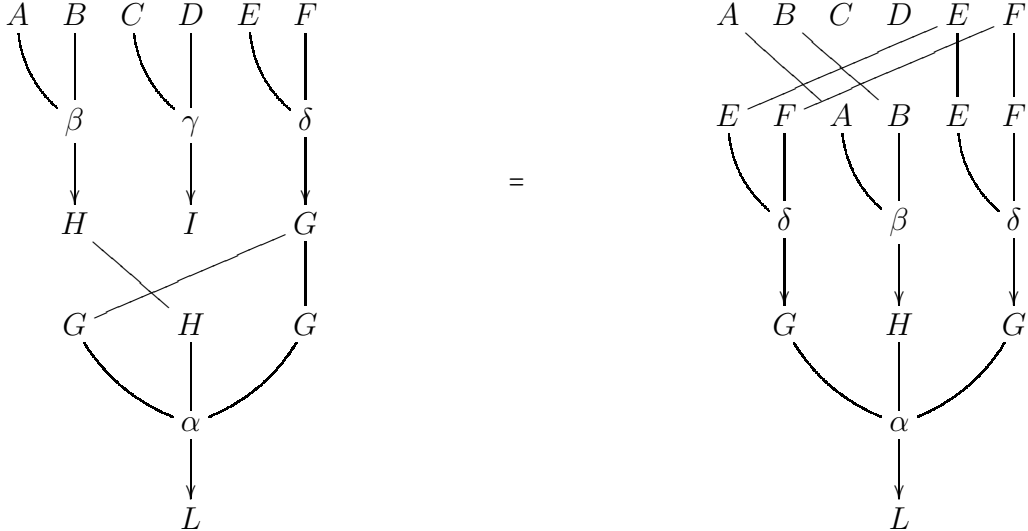
\[
\xymatrix@R=2pc@C=0.5pc{
A \ar@/_/@{-}[dr]  &   B\ar@{-}[d]   & C\ar@/_/@{-}[dr] & D\ar@{-}[d]  & E \ar@/_/@{-}[dr] &  F\ar@{-}[d]      \\
&  \beta \ar[d] && \gamma \ar[d]&& \delta \ar[d]  \\
& H   \ar@{-}[drr]         &&  I  &&   G  \ar@{-}[d]  \ar@{-}[dllll]                   \\
& G \ar@/_/@{-}[drr] &&        H   \ar@{-}[d]         & & G \ar@/^/@{-}[dll] \\
&                         &&  \al \ar[d]                           &&                           \\
&                        &&  L                                       &&
}
\xymatrix@R=2pc@C=5pc{
              &  &                       \\
              &  &                       \\
              & = & 
}
\xymatrix@R=2pc@C=0.5pc{
A \ar@{-}[drr]  &   B\ar@{-}[drr]   &  C & D &  E\ar@{-}[d] \ar@{-}[dllll] &  F\ar@{-}[d] \ar@{-}[dllll]    \\
E \ar@/_/@{-}[dr]  &   F\ar@{-}[d]   & A\ar@/_/@{-}[dr]  &  B\ar@{-}[d]  & E\ar@/_/@{-}[dr] &  F\ar@{-}[d]  &      \\
&  \delta\ar[d] && \beta\ar[d]&& \delta\ar[d]  \\
& G \ar@/_/@{-}[drr] &&        H   \ar@{-}[d]         & & G \ar@/^/@{-}[dll] \\
&                         &&  \al \ar[d]                           &&                           \\
&                        &&  L                                       &&
}
\]
As we will see in Section \ref{Cart}, to properly express this last condition one has to consider both the 
contravariant reindexing and the covariant one (see equation (\ref{FR})).
\end{enumerate}


\subsection{The category of spans in $\dM$}
\label{Span} \qq
In Proposition \ref{maincart}, we will prove the equivalence, in a cartesian multicategory $M$,
between tensor products (representability), universal products and algebraic products.
We already treated representability in Section \ref{Smc}.
While the notion itself of algebraic product requires a cartesian structure on $M$ 
(see Section \ref{Cart}), that of universal product depends only on the multicategory
structure $M:\dM \to \dPb$.
The existence of universal products is best described by the condition that the functor
$\Span(M):\Span(\dM) \to \Span(\Set_f)$, which we presently define, is an opfibration. 
In fact, $\Span(M)$ can be seen as a muticategorical version of the opposite of a fibration
(see \cite{benabou}).

Given $M:\dM\to\dPb$ in $\sMlt$ we can define the category $\Span(\dM)$,
which has the same objects as $\dM$ while an arrow is a \"span" in $\dM$, 
that is a loose arrow and a tight arrow with the same domain: 
\eq 
\xymatrix@R=1.3pc@C=1pc{
& & V \ar[ddll]_\al\ar@{-->}[ddrr]^t  & &  \\ 
& & & &  \\
X & &                  & & Y      } 
\eeq
We distinguish notationally loose and tight arrows in a diagram by naming the former with greek letters
and by drawing dashed arrows for the latter.
When appropriate, since $M_0$ is a discrete fibration, we name the tight arrows therein  
(and sometimes their domain) following the corresponding underlying mapping $f = M_0 t$:
\eq 
\xymatrix@R=1.3pc@C=1pc{
& & f^*Y \ar[ddll]_\al\ar@{-->}[ddrr]^f  & &  \\ 
& & & &  \\
X & &                  & & Y      } 
\eeq
The composition of two spans 
\eq 
\xymatrix@R=1.3pc@C=1pc{
& & V \ar[ddll]_\al\ar@{-->}[ddrr]^t  & & & & W \ar[ddll]_\beta\ar@{-->}[ddrr]^s \\ 
& & & &  \\
X & &                  & & Y  & & & & Z   } 
\eeq
is given by composing the legs in the diagram below, 
\eq 
\xymatrix@R=1.3pc@C=1pc{
& & & & U \ar[ddll]_\gamma \ar@{-->}[ddrr]^u   \\ 
& & & &  \\
& & V \ar[ddll]_\al\ar@{-->}[ddrr]^t  & & & & W \ar[ddll]_\beta\ar@{-->}[ddrr]^s \\ 
& & & &  \\
X & &                  & & Y  & & & & Z   } 
\eeq
where the cell is the unique cell obtained by lifting 
(exploiting the fact that $M_1$ is a discrete fibration)
the following pullback in $\Set_f$:
\eq 
\xymatrix@R=1.3pc@C=0.7pc{
& & U \ar[ddll] \ar[ddrr]   \\ 
& & & &  \\
MV \ar[ddrr]|{Mt}  & & & & MW \ar[ddll]|{M\beta} \\ 
& & & &  \\
& & MY  & & & &    } 
\eeq
Of course, we can treat $\Span(\Set_f)$ either as a bicategory or as a category
(by taking equivalence classes of spans).
Accordingly, we can treat $\Span(\dM)$ either as a bicategory or as a category
(by taking the suitable equivalence relation). 
We follow here this second route. 
Note that the notation $\Span(\dM)$ is slightly improper since it depends 
not only on $\dM$ but also on $M:\dM\to\dPb$. 
Anyway, we reserve $\Span(M)$ for the obvious functor $\Span(\dM) \to \Span(\Set_f)$.


\subsection{Universal products}
\label{Unipro}
The idea of a universal product in a symmetric multicategory is a very natural generalization 
of that of product in a category:
a product for $X = A_i \,\, (i\in I)$, is an object $P$ together with \"projections" $\pi_i:P\to A_i$
such that for any $Y = B_j \,\, (j\in J)$ and any family of arrows $\rho_i:Y\to A_i$, there is a unique $t:Y\to P$
with $\pi_i t = \rho_i \,\, (i\in I)$.

Here is the fibered-unbiased version of this idea.
Let $M:\dM\to\dPb$ be a symmetric multicategory.
If $X\in \dM$ and $f:pX\to J$ is a map in $\Set_f$, a {\em universal product} for $X$ along $f$ 
is an object $P\in\dM$ over $J$ with a vertical map $\pi:f^*P\to X$ 
\eq  \label{unipro1}
\xymatrix@R=1.3pc@C=1.3pc{
f^*P \ar[dd]_\pi\ar@{-->}[ddrr]^f  & &  \\ 
& &  \\
X                 & & P      } 
\xymatrix@R=2pc@C=2pc{
 \\ & \ar@{|->}[r]^M &      } \qq
\xymatrix@R=1.3pc@C=1.3pc{
pX \ar[dd]_\id\ar[ddrr]^f  & &  \\ 
& &  \\
pX                 & & J      } 
\eeq
which has the following universal property:
for any pullback $hf' = fh'$ in $\Set_f$
and any arrow $\rho: f'^*Q\to X$ over $h'$, there exists a unique cell $(t,t',f,f')$ 
over it, such that $\rho = \pi t'$.
\[
\xymatrix@R=1.3pc@C=1.3pc{
& f'^*Q \ar@{..>}[ddl]_{t'} \ar@/^/[ddddl]^(0.3)\rho \ar@{-->}[ddrr]^{f'}&& \\ \\
f^*P \ar[dd]_\pi\ar@{-->}[ddrr]^f  & & & Q \ar@{..>}[ddl]^t \\ 
& &  \\
X                 & & P      } 
\xymatrix@R=2pc@C=2pc{
 \\ \\ & \ar@{|->}[r]^M &      } \qq
\xymatrix@R=1.3pc@C=1.3pc{
& L \ar[ddl]_{h'} \ar@/^/[ddddl]^(0.3){h'} \ar[ddrr]^{f'}&& \\ \\
pX \ar[dd]_\id\ar[ddrr]^f  & & \pb & K \ar[ddl]^h \\ 
& &  \\
pX                 & & J      } 
\]

\begin{remark}
It is easy to see that the above sketched notion of universal product in $M$ arises as a special case 
of the present one, namely by taking the product of $X = A_i \,\, (i\in I)$ along $f:I\to 1$.
\end{remark}

\begin{prop}
\label{unipro3}
The multicategory $M:\dM\to\dPb$ has universal products if and only if 
$\Span(M):\Span(\dM) \to \Span(\Set_f)$ is an opfibration.
\end{prop}
\pf
It is easy to see that the condition above for (\ref{unipro1}) being an universal product, 
amounts to the condition of being an opcartesian arrow for $\Span(M)$ over $f$, 
where the unique lifting condition is restricted to spans in $\Set_f$ whose right leg is an identity.
From this, the unique lifting condition for general arrows follows immediately. 
\epf


\subsection{The multicategory of enhanced spans}
\label{Espan} \qq
In order to define the monad $(-)^\cart$ we need to define the
category $\Espan(\C)$ of \"enhanced spans" in a category $\C$ with pullbacks.
We use this provisional terminology, since we don't now if this notion has been 
already considered and named in the literature. 

The objects of $\Espan(\C)$ are the objects of $\C$, while the arrows $A\to B$
are the commutative triangles 
\eq 
\xymatrix@R=1.3pc@C=1pc{
& & C \ar[ddll]_b \ar[ddrr]^a  & &  \\ 
& & & &  \\
B  & & & & A \ar[llll]^f     } 
\eeq
and composition is given by pullbacks:
\eq 
\xymatrix@R=1.3pc@C=1pc{
& & & & F \ar[ddll]_r \ar[ddrr]^s   \\ 
& & & &  \\
& & D \ar[ddll]_c\ar[ddrr]^b  & & \pb & & E \ar[ddll]_{b'}\ar[ddrr]^a \\ 
& & & &  \\
C & &                  & & B \ar[llll]^g & & & & A \ar[llll]^f  } 
\eeq
It is easy to see that $\Espan(\C)$ is the loose part of a (non-strict) double category 
$\dEspan(\C)$ whose tight part is $\C$ itself and whose cells are given by pullbacks
(the three rectangles): 
\eq 
\xymatrix@R=1.3pc@C=1pc{
& & C' \ar[dd] \ar[ddll] \ar[ddrr]  & &  \\
& & & &  \\
B' \ar[dd] & & C \ar[ddll] \ar[ddrr]  & & A' \ar[dd]\ar@/^1pc/[llll] \\ 
& & & &  \\
B  & & & & A \ar@/^/[llll]    } 
\eeq
It is also easy to see that, for $\C = \Set_f$, the obvious functor $N:\dEspan \to \dPb$
\eq 
\xymatrix@R=1.3pc@C=1pc{
& & C' \ar[dd] \ar[ddll] \ar[ddrr]  & &  \\
& & & &  \\
B' \ar[dd] & & C \ar[ddll] \ar[ddrr]  & & A' \ar[dd]\ar@/^1pc/[llll] \\ 
& & & &  \\
B  & & & & A \ar@/^/[llll]    } 
\xymatrix@R=2pc@C=2pc{
 \\ \\ & \ar@{|->}[r]^N &      } \qq
 \xymatrix@R=1.3pc@C=1pc{
\\ \\
B' \ar[dd] & & & & A' \ar[dd]\ar@/^1pc/[llll] \\ 
& & & &  \\
B  & & & & A \ar@/^/[llll]    } 
\eeq
is a (non-strict) unbiased symmetric multicategory. 
Indeed, it is equivalent to $\dN\t$, where $\dN$ is the monoid of natural numbers under multiplication.

Now, for any symmetric multicategory $M:\dM \to \dPb$, 
we can define a double category $\dEspan(\dM)$ whose tight arrows are those of $\dM$ and whose
loose arrows are spans in $\dM$ over enhanced spans in $\Set_f$.
There is an obvious functor $\dEspan(\dM) \to \dEspan$ which, composed with $N:\dEspan \to \dPb$
gives a symmetric multicategory $M^\cart: \dEspan(\dM) \to \dPb$.


\subsection{The monad $(-)^\cart$}
\label{Moncart} \qq
Thus, for any symmetric multicategory $M:\dM \to \dPb$ we have defined another 
symmetric multicategory $M^\cart: \dEspan(\dM) \to \dPb$, and this construction extends 
to a monad $(-)^\cart$ on $\sMlt$.
Indeed, since the loose arrows of $M^\cart$ over $b:I\to J$ are the spans in $\dM$
over commutative triangles:
\eq 
\xymatrix@R=1.3pc@C=1pc{
& & V \ar[ddll]_\al \ar@{-->}[ddrr]^f  & &  \\ 
& & & &  \\
Y & & & & X     } 
\xymatrix@R=2pc@C=2pc{
 \\ & \ar@{|->}[r] &      } \qq
\xymatrix@R=1.3pc@C=1pc{
& & K \ar[ddll]_a \ar[ddrr]^f  & &  \\ 
& & & &  \\
J  & & & & I \ar[llll]^b     } 
\eeq
the loose arrows of $(M^\cart)^\cart$ over $c: L\to J$ have the form:
\eq 
\xymatrix@R=1.3pc@C=1pc{
& & V \ar[ddll]_\al \ar@{-->}[dr]^f  & &  \\ 
& & & X \ar@{-->}[dr]^g &  \\
Y & & & & Z     } 
\xymatrix@R=2pc@C=2pc{
 \\ & \ar@{|->}[r] &      } \qq
\xymatrix@R=1.3pc@C=1pc{
& & K \ar[ddll]_a \ar[dr]^f   & &  \\ 
& & & I \ar[dr]^g \ar[dlll]_b &  \\
J & & & & L \ar[llll]^c    } 
\eeq
It is then clear that composing the right legs of the spans
\eq 
\xymatrix@R=1.3pc@C=1pc{
& & V \ar[ddll]_\al \ar@{-->}[dr]^f & &  \\ 
& & & X \ar@{-->}[dr]^g &  \\
Y & & & & Z     } 
\xymatrix@R=2pc@C=2pc{
 \\ & \ar@{|->}[r]^\mu &      } \qq
\xymatrix@R=1.3pc@C=1pc{
& & V \ar[ddll]_\al \ar@{-->}[ddrr]^{g\circ f} & &  \\ 
& & & &  \\
Y & & & & Z     } 
\eeq
gives a morphism $\mu:(M^\cart)^\cart \to M^\cart$
which, along with the obvious $\eta:M \to M^\cart$
\eq 
\xymatrix@R=1.3pc@C=1pc{
& & X \ar[ddll]_\al  & &  \\ 
& & & &  \\
Y & & & &      } 
\xymatrix@R=2pc@C=2pc{
 \\ & \ar@{|->}[r]^\eta &      } \qq
\xymatrix@R=1.3pc@C=1pc{
& & V \ar[ddll]_\al \ar@{-->}[ddrr]^{\id} & &  \\ 
& & & &  \\
Y & & & & X     } 
\eeq
defines the desired monad on $\sMlt$


\subsection{Cartesian multicategories and the covariant reindexing}
\label{Cart} \qq
We define the category $\sMlt^\cart$ of {\em cartesian multicategories} as the 
category of algebras for the monad $(-)^\cart:\sMlt\to\sMlt$.
Thus a cartesian structure $\Gamma$ on $M:\dM \to \dPb$ consists of a morphism
$\Gamma:M^\cart \to M$ in $\sMlt$ satisfying the usual conditions for monad algebras.

We presently make explicit these conditions and show how they reduce to the usual
defining axioms for a cartesian multicategory sketched in Section \ref{Classic}.

If $(M,\Gamma)$ is a cartesian multicategory, for any enhanced span $(\al,f)$ in $\dM$ 
(that is, for any loose arrow in $\dEspan(\dM)$)
we have a loose arrow $\beta = \Gamma_l(\al,f)$ in $\dM$, which we interpret as
the {\em covariant reindexing} of $\al$ along $f$ 
(over the commutative triangle $hf = M\al$ in $\Set_f$):
\eq \label{cov}
\xymatrix@R=2pc@C=2pc{
X \ar[ddr]_\al\ar@{-->}[rr]^f  & & Y \ar@{..>}[ddl]^\beta \\ 
& &  & \ar@{|->}[r]^M & \\
             & Z &      }    \qq
\xymatrix@R=1.8pc@C=1.8pc{
MX \ar[ddr]_{M\al}\ar[rr]^f  & & I  \ar[ddl]^h \\ 
&& \\
             & MZ &      }
\eeq
We will often denote such a $\beta = \Gamma_l(\al,f)$ also by $f_!\al $:
\eq
f_!\al = \Gamma_l(\al,f)
\eeq
Note however that both notations are somewhat imprecise, 
since $\beta$ depends also on the underlying enhanced span.

In order to emphasize, in a diagram, that $\beta = \Gamma_l(\al,f)$,
we put the symbol $\#$ inside the corresponding triangle:
\eq
\xymatrix@R=2pc@C=2pc{
X \ar@/_/[ddr]_\al\ar@{-->}[rr]^f  & & Y \ar@/^/[ddl]^{\beta = f_!\al} \\ 
& \# &  \\
& Z &      }   
\eeq
Similarly, we emphasize cells in $\dM$ by putting the symbol $@$ inside the corresponding square:
\eq
\xymatrix@R=2pc@C=2pc{
X' \ar[dd]_{\al=f^*\beta}\ar@{-->}[rr]^{f'}  & & Y' \ar[dd]^\beta \\ 
&  @ & \\
X \ar@{-->}[rr]^f  & & Y  
}   
\eeq
In this case, we write $\al = f^*\beta$. As for $f_!$, this notation is somewhat imprecise,
since $\al$ depends also on the underlying pullback.
\begin{prop}
\label{funct}
The covariant reindexing is functorial: 
\[ 
f_!(g_! \al) = (fg)_! \al 
\]
\end{prop}
\pf
It is a direct consequence of the general definition of algebras for a monad 
and of the definition of the monad $(-)^\cart$.
\epf
Functoriality can be rephrased by saying that covariant reindexings can be pasted horizontally: 
\eq
\xymatrix@R=2pc@C=1pc{
X \ar@/_1pc/[ddrr]_\al\ar@{-->}[rr]^f  & & Y \ar@{-->}[rr]^g \ar[dd]^\beta & & Z \ar@/^1pc/[ddll]^\gamma \\ 
& \# & & \# \\
& & W &      }  
\xymatrix@R=2pc@C=2pc{
\\ & \ar@{=>}[r] &      } \qq
 \xymatrix@R=2pc@C=1pc{
X \ar@/_/[ddrr]_\al\ar@{-->}[rrrr]^{gf}  & & & & Z \ar@/^/[ddll]^\gamma \\ 
& & \# \\
& & W &      }  
\eeq
\begin{remark}
The explicit proof of Proposition \ref{funct} is very similar to the proof that the algebras for the 
monad $M\otimes-$, where $M$ is a monoid in a monoidal category, are the actions of $M$.
Indeed, the monad $(-)^\cart$ itself is similar to the monad $M\times-$ (or $-\times M$) on sets:
so as the latter takes a set $X$ to the set of pairs $(x,m)$, 
the former takes the set of loose arrows of $\dM$ to the set of spans $(\al,f)$, 
obtained by pairing them with tight arrows.
Only that, in this case, the set is in fact a category, so that also the spans $(\al,f)$ should
be the arrows of a category. 
This, as explained above, is achieved by exploiting the double category structure of $\dM$,
which allows to define $\dEspan(\dM)$.
Even if this construction suggests the presence of two monads related by 
a distributive law, I don't know yet how to make explicit this intuition.
\end{remark}


\begin{prop}
The functoriality of the loose part $\Gamma_l$ of $\Gamma$ implies to the following conditions: 
\begin{enumerate}
\item
Pasting a covariant and contravariant reindexing gives a covariant reindexing.
\item
Tailing a covariant reindexing with a loose arrow, gives again a covariant reindexing.
\end{enumerate}
\end{prop}
\pf
Recall that the composition $(\beta,g)\circ(\al,f)$ in $M^\cart$ is obtained as follows:
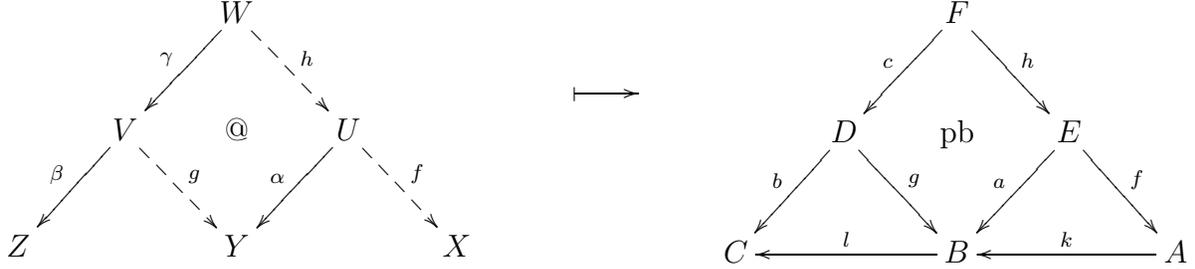
\eq 
\nonumber
\xymatrix@R=1pc@C=0.8pc{
& & & & W \ar[ddll]_\gamma \ar@{-->}[ddrr]^h   \\ 
& & & &  \\
& & V \ar[ddll]_\beta\ar@{-->}[ddrr]^g  & & @ & & U \ar[ddll]_\al\ar@{-->}[ddrr]^f \\ 
& & & &  \\
Z & &                  & & Y  & & & & X   } 
\xymatrix@R=2pc@C=2pc{
\\  & \ar@{|->}[r] &      } \qq
\xymatrix@R=1pc@C=0.8pc{
& & & & F \ar[ddll]_c \ar[ddrr]^h   \\ 
& & & &  \\
& & D \ar[ddll]_b\ar[ddrr]^g  & & \pb & & E \ar[ddll]_a\ar[ddrr]^f \\ 
& & & &  \\
C & &                  & & B \ar[llll]_l & & & & A \ar[llll]_k  } 
\eeq

Consider the particular case $f = \id_X$, that is the composition
$(\beta,g)\circ(\al,\id)$, and say that $\Gamma_l(\beta,g) = \delta$. 
Since $\Gamma_l(\al,\id) = \Gamma_l(\eta\al) = \al$, 
by the unit condition for the algebra $M$, we have
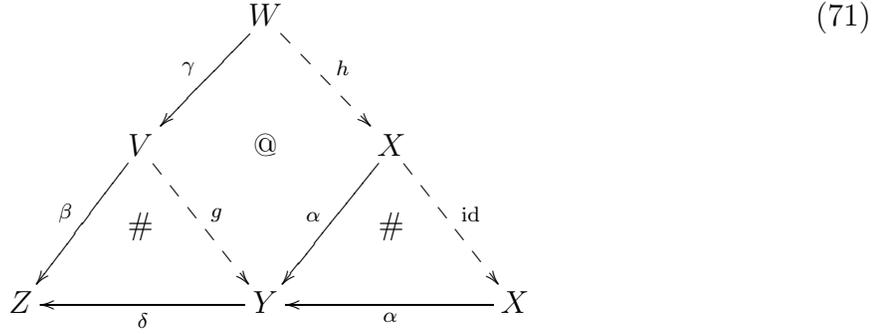
\eq 
\xymatrix@R=1.2pc@C=1pc{
& & & & W \ar[ddll]_\gamma \ar@{-->}[ddrr]^h   \\ 
& & & &  \\
& & V \ar[ddll]_\beta\ar@{-->}[ddrr]^g  & & @ & & X \ar[ddll]_\al\ar@{-->}[ddrr]^\id \\ 
& & \#  & & & & \# \\
Z & &    & & Y  \ar[llll]^\delta  & & & & X  \ar[llll]^\al   } 
\eeq

Since $\Gamma_l$ preserves composition, we have 
\[
\Gamma_l(\beta\gamma,h) = \Gamma_l((\beta,g)\circ(\al,\id)) =
\Gamma_l(\beta,g)\circ \Gamma_l(\al,\id) = \delta\al 
\]
that is the following \"pasting" property holds:
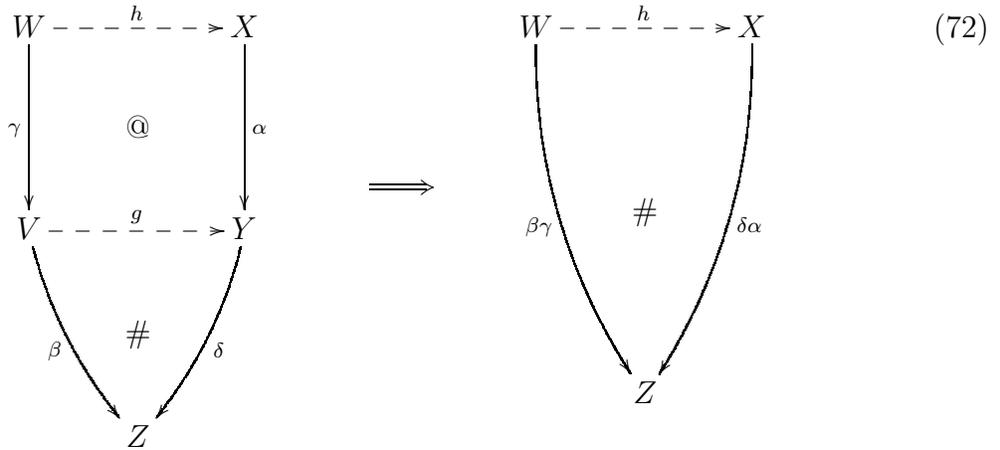
\eq
\label{FR}
\xymatrix@R=2pc@C=2pc{
W \ar[dd]_\gamma\ar@{-->}[rr]^h  & & X \ar[dd]^\al \\ 
&  @ & \\
V \ar@/_/[ddr]_\beta\ar@{-->}[rr]^g  & & Y \ar@/^/[ddl]^\delta \\ 
& \# &  \\
& Z &      }  
\xymatrix@R=2pc@C=2pc{
\\ \\ & \ar@{=>}[r] &      } \qq
\xymatrix@R=2pc@C=2pc{
W \ar@/_1pc/[ddddr]_{\beta\gamma}\ar@{-->}[rr]^h  & & X \ar@/^1pc/[ddddl]^{\delta\al} \\ 
&  & \\
&  \# & \\
&  &  \\
& Z &      }  
\eeq

Consider now the particular case $g = \id_Y$, that is the composition
$(\beta,\id)\circ(\al,f)$, and say that $\Gamma_l(\al,f) = \delta$:
\eq 
\xymatrix@R=1.2pc@C=1pc{
& & & & W \ar[ddll]_\al \ar@{-->}[ddrr]^\id   \\ 
& & & &  \\
& & Y \ar[ddll]_\beta\ar@{-->}[ddrr]^\id  & & @ & & U \ar[ddll]_\al\ar@{-->}[ddrr]^f \\ 
& & \#  & & & & \# \\
Z & &    & & Y  \ar[llll]^\beta  & & & & X  \ar[llll]^\delta   } 
\eeq
Since $\Gamma_l$ preserves composition we have 
\[
\Gamma_l(\beta\alpha,f) = \Gamma_l((\beta,\id)\circ(\al,f)) =
\Gamma_l(\beta,\id)\circ \Gamma_l(\al,f) = \beta\delta 
\]
that is the following \"tailing" property for covariant reindexing holds, 
for any loose arrow $\beta:Y\to Z$ in $\dM$:
\eq
\label{TAIL}
\xymatrix@R=2pc@C=2pc{
U \ar@/_/[ddr]_\al\ar@{-->}[rr]^f  & & X \ar@/^/[ddl]^\delta \\ 
& \# &  \\
& Y  \ar[dd]_\beta &  \\
&  & \\  
& Z & }  
\xymatrix@R=2pc@C=2pc{
\\ \\ & \ar@{=>}[r] &      } \qq
\xymatrix@R=2pc@C=2pc{
U \ar@/_1pc/[ddddr]_{\beta\al}\ar@{-->}[rr]^f  & & X \ar@/^1pc/[ddddl]^{\beta\delta} \\ 
&  & \\
&  \# & \\
&  &  \\
& Z &      }  
\eeq

\epf

We refer to (\ref{FR}) as the \"Frobenius law". Indeed it may be written as
\[ h_! (\beta(f^*\al)) = (g_!\beta)\al \]
which is indeed very similar to the classical Frobenius reciprocity law.

On the other hand, the functoriality of $\Gamma$ on cells says that covariant reindexings 
are stable with respect to the contravariant ones:
if the front and back rectangles below are cells in $\dM$ and the triangle on the right is a covariant reindexing,
so it is also the triangle on the left.
\eq
\label{BC}
\xymatrix@R=1.5pc@C=1.5pc{
&& Y' \ar[dddl] \ar@{-->}[rrrr] &&&& Y\ar[dddl] \\
X' \ar[ddr]\ar@{-->}[rru]^{f'}  \ar@{-->}[rrrr] &&&& X\ar[ddr]_a\ar@{-->}[rru]^f  \\ 
&&  \\
& Z' \ar@{-->}[rrrr]_g &&&& Z &      }   
\eeq

We refer to (\ref{BC}) as the \"Beck-Chevalley law". Indeed it may be written as
\[ g^*(f_!a) = f'_!(g^*a) \]
which is indeed a double categorical form of the usual Beck-Chevalley law for fibrations.

Summarizing, we can give the following equivalent 
\begin{definition}
\label{defcart2}
A {\em cartesian multicategory} is an unbiased symmetric multicategory
with a covariant functorial reindexing (along tight arrows and over commutative triangles in $\Set_f$) 
which is related to the contravariant one by the Frobenius and the Beck-Chevalley laws
and such that the tailing condition holds.
\end{definition}
This definition is also equivalent to the classical definition as given for instance in \cite{pisani}.
Indeed, the conditions in \ref{defcart2} are the fibrational (unbiased) form 
of the classical conditions.
In fact, it is easy to see that the tailing and the Frobenius conditions correspond respectively 
to the compatibility of reindexing with composition from above and from below.
On the other hand, the Beck-Chevalley Law assures that the reindexing of a family of maps is the 
family of the reindexed maps. 

\begin{remark}
\label{rmkfrob}
It is clear, in the author's opinion, that the Frobenius Law is the right way to express
compatibility from below.
The passage from symmetric to cartesian multicategories is not merely a matter of extending
a covariant reindexing from bijective mapping to arbitrary mapping as is sometimes suggested:
both the covariant and the contravariant reindexing are involved in the notion of cartesian multicategory.
\end{remark}

\begin{remark}
In \cite{pisani2} we adopted essentially Definition \ref{defcart2},
which allowed to prove, in the fibrational setting, Proposition \ref{maincart} 
(as expected from any sensible definition).
Now, the fact that it is equivalent, for a multicategory $M$, to being the algebra of a rather natural monad
enforces the belief that we have so captured a suitable notion of cartesian multicategory.
\end{remark}


\subsection{Coherence of the two reindexing}
One may wonder whether the two reindexing are coherent, when both are possible.
This is in fact guaranteed by the Beck-Chevalley law, as shown in the following proposition.
Given a cartesian fibered multicategory $\dM$, a commutative triangle in $\Set_f$ whose top side is 
an {\em isomorphism}, and a lifting $\al:X\to Y$ of its left side in $\dM$,
we can get not only the covariant reindexing $f_!\al$, but also the contravariant
reindexing $g^*\al$ along the inverse map $g = f^{-1}$:
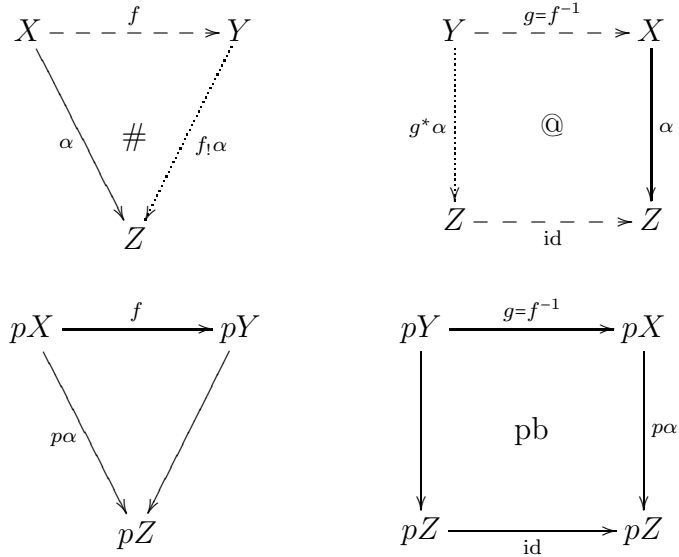
\[
\xymatrix@R=2pc@C=2pc{
X \ar[ddr]_\al\ar@{-->}[rr]^f  & & Y \ar@{..>}[ddl]^{f_!\al} \\ 
& \# & & \\
& Z &      }    \qq
\xymatrix@R=1.8pc@C=1.8pc{
Y \ar@{..>}[dd]_{g^*\al}\ar@{-->}[rr]^{g=f^{-1}}  & & X \ar[dd]^\al \\ 
& @ & \\
Z \ar@{-->}[rr]_\id                 & & Z      }             
\]
\[
\xymatrix@R=2.3pc@C=1.5pc{
pX \ar[ddr]_{p\al}\ar[rr]^f  & & pY  \ar[ddl] \\ 
&& \\
             & pZ &      } \qq\qq
\xymatrix@R=1.8pc@C=1.8pc{
pY \ar[dd]\ar[rr]^{g=f^{-1}}  & & pX \ar[dd]^{p\al} \\
& {\rm pb} & \\
pZ  \ar[rr]_\id                 & & pZ       }             
\]

\begin{prop}
In the above situation, the two reindexing coincide: $f_!\al = g^*\al$.
\end{prop}
\pf
Consider the following prism, where $g = f^{-1}$ (so that the top square commutes)
and where we assume that the right triangle is a covariant reindexing and 
that the front and back rectangles are cells in $\dM$.
Then, by the Beck-Chevalley condition, the left triangle is also a covariant reindexing, so that $t = g^*\al$.
But since the back rectangle is a cell, also $t = f_!\al$.
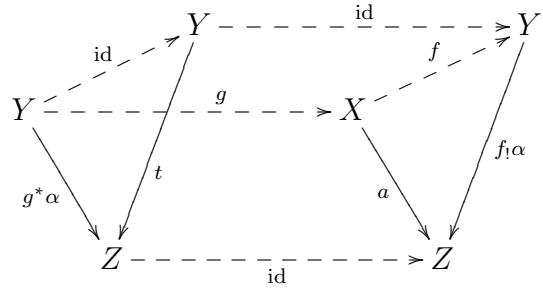
\[
\xymatrix@R=1.5pc@C=1.5pc{
&& Y \ar[dddl]^(.6)t \ar@{-->}[rrrr]^\id &&&& Y\ar[dddl]^{f_!\al} \\
Y \ar[ddr]_{g^*\al}\ar@{-->}[rru]^\id  \ar@{-->}[rrrr]^(.6)g &&&& X\ar[ddr]_a\ar@{-->}[rru]^(.5)f  \\ 
&&  \\
& Z \ar@{-->}[rrrr]_\id &&&& Z &      }   
\]
\epf


\subsection{Algebraic products}
The idea of a cartesian multicategory is that it is a \"virtual finite product category": 
if it is representable, the tensor product is cartesian. 
This is made precise in Theorem \ref{maincart}, where it is shown that 
a cartesian multicategory is representable if and only if it has universal products
(as defined in Section \ref{Unipro}). 
Indeed, both conditions are equivalent to the existence of \"algebraic products", that are defined
equationally in cartesian multicategories and are thus \"absolute" 
(that is, are preserved by any cartesian functor).

If $X\in\dM$ and $f:pX\to J$ is a map in $\Set_f$, an {\em algebraic product} for $X$ along $f$ 
is an object $P\in\dM$ over $J$ along with a vertical map $\pi:f^*P\to X$ and a map $u:X\to P$ over $f$
\eq \label{algpr1}
\xymatrix@R=1.3pc@C=1.3pc{
f^*P \ar[dd]_\pi\ar@{-->}[ddrr]^f  & &  \\ 
& & \\
X    \ar[rr]_u            & & P      } \q
\xymatrix@R=1.3pc@C=1.3pc{
& &  \\ 
&\ar@{|->}[r] & \\
& &      } \qq
\xymatrix@R=1.3pc@C=1.6pc{
pX \ar[dd]_\id\ar[ddrr]^f  & &  \\ 
& &  \\
pX  \ar[rr]_f               & & J      } 
\eeq
such that the following are both covariant reindexing triangles:
\eq. \label{algpr2}
\xymatrix@R=2pc@C=0.5pc{
f^*P \ar[d]_\pi\ar@{-->}[rr]^f  & & P \ar@/^1pc/[ddl]^\id \\ 
X \ar@/_/[dr]_u & \# & & \\
& P &      }    \qq \qq \qq
\xymatrix@R=2pc@C=0.5pc{
X \ar@/_1pc/[ddr]_\id\ar@{-->}[rr]^\Delta  & & h^*X \ar[d]^{f^*u} \\ 
& \# & f^*P \ar@/^/[dl]^\pi & \\
& X &      }    
\eeq
where the map $\Delta$ in the right hand triangle is a tight arrow over 
the diagonal of the pullback of $f$ along itself in $\Set_f$,
so that $h\Delta = \id$ also in $\dM_0$:
\[
\xymatrix@R=1.2pc@C=1pc{
X  \ar@{-->}[dr]^\Delta \ar@{-->}@/^0.9pc/[drrr]^\id &&& \\
& h^*X \ar[dd]_{f^*u}\ar@{-->}[rr]^h  & & X \ar[dd]^u \\
& & @ & \\
& f^*P  \ar@{-->}[rr]_f                 & & P       }     \qq
\xymatrix@R=1.3pc@C=1.3pc{
& &  \\ \\
&\ar@{~>}[r] & \\
& &      } \qq
\xymatrix@R=1pc@C=1pc{
pX  \ar[dr]^\Delta  \ar@/_0.7pc/[dddr]_\id \ar@/^0.7pc/[drrr]^\id &&& \\
& K \ar[dd]^l\ar[rr]^h  & & pX \ar[dd]^f \\
& & \pb & \\
& pX  \ar[rr]_f                 & & J       }             
\]

Now we come to the main result of this section.

\begin{theorem}
\label{maincart}
For a cartesian multicategory $M:\dM\to\dPb$, the following are equivalent:
\begin{enumerate}
\item
$M$ has algebraic products. {\q \em (AP)}
\item
$M$ has universal products. {\q \em (UP)}
\item
$M$ is representable. {\q \em (R)}
\end{enumerate}
\end{theorem}
\pf 
Recall that we use to mark the covariant reindexing triangles with the symbol $\#$ inside 
and the contravariant reindexing squares (the cells of $\dM$) with the symbol $@$ inside.
In the present proof, we refer to them respectively as the \"special triangles" and the \"special squares".
The wavy arrows relate a diagram in $\dM$ with the underlying diagram in $\dPb$.

Let us show that (AP) implies (UP). Suppose then that $(\pi,u)$ is an algebraic product
for $X$ along $f$ as in (\ref{algpr1}). 
We want to show that $\pi$ is a universal product for $X$ along $f$ as in (\ref{unipro1}):
given a pullback $hf' = fh'$ in $\Set_f$ and a tight lifting of its top side $f'$,
any loose arrow $\rho: f'^*Q\to X$ with $p\rho = h'$, 
should factor uniquely as $\rho = \pi(f^*t)$ for a unique special square: 
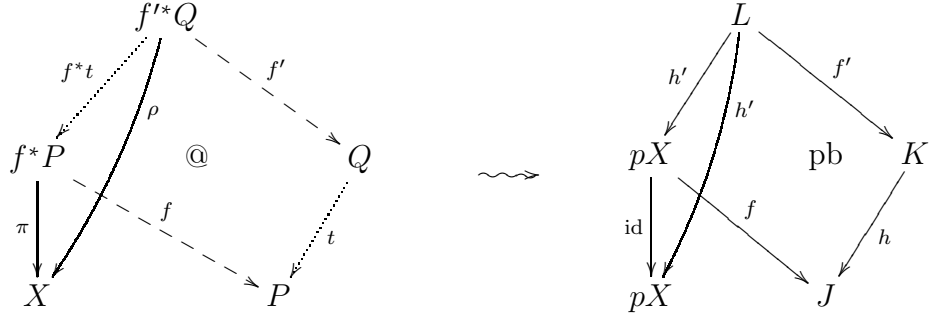
\[
\xymatrix@R=1.3pc@C=1.3pc{
& f'^*Q \ar@{..>}[ddl]_{f^*t} \ar@/^/[ddddl]^(0.3)\rho \ar@{-->}[ddrr]^{f'}&& \\ \\
f^*P \ar[dd]_\pi\ar@{-->}[ddrr]^f  &  \qq@ & & Q \ar@{..>}[ddl]^t \\ 
& &  \\
X                 & & P      } 
\xymatrix@R=2pc@C=2pc{
 \\ \\ & \ar@{~>}[r] &      } \qq
\xymatrix@R=1.3pc@C=1.3pc{
& L \ar[ddl]_{h'} \ar@/^/[ddddl]^(0.3){h'} \ar[ddrr]^{f'}&& \\ \\
pX \ar[dd]_\id\ar[ddrr]^f  & & \pb & K \ar[ddl]^h \\ 
& &  \\
pX                 & & J      } 
\]

Supposing that such a factorization exists, by the Frobenius law applied to
the first of (\ref{algpr2}), we see that $t = f'_!(u\rho)$:
\[
\xymatrix@R=2pc@C=1.5pc{
f'^*Q \ar[dd]_{f^*t}\ar@{-->}[rr]^{f'}  & & Q \ar[dd]^t \\ 
& @ & \\
f^*P \ar[d]_\pi\ar@{-->}[rr]^f  & & P \ar@/^0.8pc/[ddl]^\id \\ 
X \ar@/_/[dr]_u & \# & & \\
& P &      }    \qq \qq \qq
\xymatrix@R=2pc@C=1.5pc{
f'^*Q \ar[dd]_\rho\ar@{-->}[rr]^{f'}  & & Q \ar@/^1.3pc/[ddddl]^t \\ \\
X \ar@/_/[ddr]_u & \# & & \\ \\
& P &      } 
 \]
This proves unicity. To show that such a $t$ gives indeed the desired factorization,
consider the diagram below:
\[
\xymatrix@R=1pc@C=1pc{
&&& f'^*Q\ar@{-->}[rr]^{f'}\ar@/^1pc/[dddddl]|(.4){f^*t} & & Q  \ar@/^1pc/[dddddl]^t  \\
f'^*Q\ar@{-->}@/^0.8pc/[urrr]^\id\ar@{-->}[rr]^\Delta \ar[dd]_\rho && 
W \ar@{-->}[rr]\ar@{-->}[ur]\ar[dd]_{h^*\rho} && f'^*Q \ar[dd]_\rho\ar@{-->}[ur]  \\ 
&&&& & & & \\
X\ar@{-->}[rr]^\Delta\ar@/_1pc/[rrdddd]_\id && h^*X \ar@{-->}[rr]^h\ar[dd]_{f^*u} && X \ar[dd]_u & & & \\ \\
&& f^*P \ar@{-->}[rr]^f\ar[dd]_\pi && P &&&    \\ \\  
&& X }
\]
The right-hand triangle (with $t$ as a side) is special by hypothesis, 
so by Beck-Chevalley also the triangle with $f^*t$ as a side is special.
The lower left-hand triangle is special by the second of (\ref{algpr2}) so that,
by Frobenius, it is special also its pasting with the left-hand special square.
Now, by composition, we get a special triangle with an identity top side,
so that $\rho = \pi(f^*t)$, as desired.

Next, we prove that (AP) implies (R). First let us show that
in an algebraic product $(\pi,u)$ along $f:pX\to J$, $u$ is opcartesian in $\dM_l$ over $f$.
So let $v:X\to Q$ and let $g: J\to pQ$ be such that $gf = pv$;
then there should be a unique map $t:P\to Q$ over $g$, such that $tu = v$:
\[
\xymatrix@R=1.3pc@C=1.3pc{
X  \ar[ddrr]_v  \ar[rr]^u    & & P \ar@{..>}[dd]^t  \\
& & \\
& & Q  } \q
\xymatrix@R=1.3pc@C=1.3pc{
& &  \\ 
&\ar@{~>}[r] & \\
& &      } \qq
\xymatrix@R=1.3pc@C=1.3pc{
pX  \ar[ddrr]_{pv}  \ar[rr]^f    & & J \ar[dd]^g  \\
& & \\
& & pQ  }
\]
Supposing that such a factorization exists, by the condition (1) in the definition
of cartesian multicategory applied to the first of (\ref{algpr2}), we see that $t = f_!(v\pi)$:
\[
\xymatrix@R=2pc@C=1.5pc{
f^*P \ar[d]_\pi\ar@{-->}[rr]^f  & & P \ar@/^0.8pc/[ddl]^\id \\ 
X \ar@/_/[dr]_u & \# & & \\
& P \ar[d]_t &   \\  
& Q &   }    \qq \qq \qq
\xymatrix@R=1pc@C=1.5pc{
f^*P \ar[dd]_\pi\ar@{-->}[rr]^f  & & P \ar@/^1pc/[ddddl]^t \\ \\
X \ar@/_/[ddr]_v & \# & & \\ \\
& P &      } 
 \]
This proves unicity. To show that such a $t$ gives indeed the desired factorization,
consider the diagram below:
\[
\xymatrix@R=1pc@C=1.5pc{
X \ar@{-->}[rr]^\Delta \ar@/_1pc/[ddddrr]_\id \ar@/^2pc/@{-->}[rrrr]^\id 
&& h^*X \ar[dd]_{f^*u}\ar@{-->}[rr]^h  & & X \ar[dd]^u \\ 
& \# & &@ & \\
&& f^*P \ar[dd]_\pi\ar@{-->}[rr]^f  & & P \ar@/^1pc/[ddddl]^t \\ \\
&& X \ar@/_/[ddr]_v & \# & & \\ \\
&& & P &      } 
 \]
The outer triangle is special (by pasting and compositions) and has an identity top side,
so that $v = tu$, as desired.
To show that $u$ is {\em stably} opcartesian, note that the notion of algebraic product is 
stable with respect to reindexing: 
if $(\pi,u)$ is an algebraic product along $f:pX\to J$ and $l:L\to J$ is a map in $\Set_f$,
then reindexing the (\ref{algpr1}) along $l$ one gets an algebraic product $(\pi',u')$
along $f'$ (the pullback of $f$ along $l$).
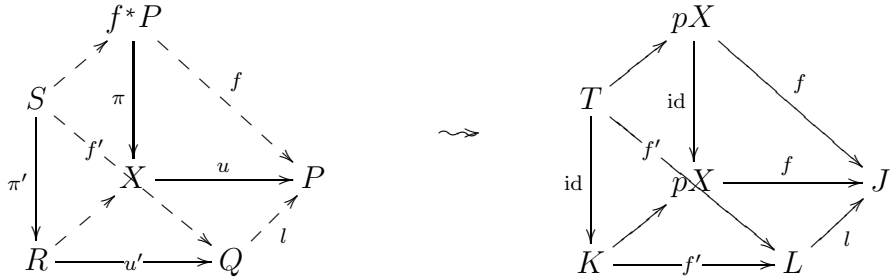
\[
\xymatrix@R=1.3pc@C=1.3pc{
& f^*P \ar[dd]_\pi\ar@{-->}[ddrr]^f  & &  \\ 
S \ar@{-->}[ur]\ar[dd]_{\pi'} \ar@{-->}[ddrr]|(.3){f'} & & \\
& X    \ar[rr]^u            & & P   \\
R \ar@{-->}[ur]\ar[rr]|{u'} && Q \ar@{-->}[ur]_l    } \q
\xymatrix@R=1.3pc@C=1.3pc{
& &  \\  \\
&\ar@{~>}[r] & \\
& &      } \qq
\xymatrix@R=1.3pc@C=1.6pc{
& pX \ar[dd]_\id\ar[ddrr]^f  & &  \\ 
T \ar[ur]\ar[dd]_\id\ar[ddrr]|(.3){f'} & &  \\
& pX  \ar[rr]^f               & & J  \\
K \ar[ur]\ar[rr]|{f'} && L \ar[ur]_l   } 
\]
Indeed, since (contravariant) reindexing preserves identity, composition, special squares
and special triangles, it preserves the \"equations" (\ref{algpr2}) as well. 

Next, we prove that (UP) implies (AP). Suppose then that $\pi:f^*P\to X$
is a universal product for $X$ along $f$:
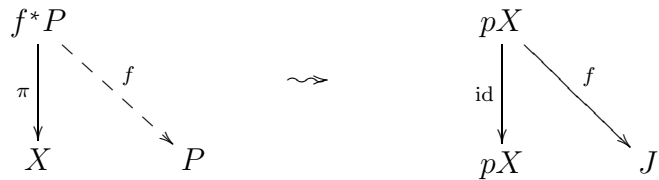
\[
\xymatrix@R=1.3pc@C=1.3pc{
f^*P \ar[dd]_\pi\ar@{-->}[ddrr]^f  & &  \\ 
& &  \\
X                 & & P      } 
\xymatrix@R=1.3pc@C=1.3pc{
& &  \\ 
&\ar@{~>}[r] & \\
& &      } \qq\qq
\xymatrix@R=1.3pc@C=1.3pc{
pX \ar[dd]_\id\ar[ddrr]^f  & &  \\ 
& &  \\
pX                 & & J      } 
\]
We show that $\pi$ is part of an algebraic product $(\pi,u)$. 
The map $u:X\to P$ is obtained by exploiting the universal property of $\pi$
(in fact, the \"existence" part).
It is the (unique) map such that $\pi(f^*u) = \Delta_!\id$, that is such that
both the square and the triangle in the diagram below are special:
\[
\xymatrix@R=4pc@C=1.5pc{
X \ar@/_1pc/[ddr]_\id \ar@{-->}[rr]^\Delta  & & h^*X \ar@{..>}[d]^{f^*u}\ar@{-->}[rr]^h && X \ar@{..>}[d]^u \\ 
&  & f^*P \ar@/^.8pc/[dl]^\pi \ar@{-->}[rr]_f & & P \\
& X &      }    
\]
To show that also the first of equations (\ref{algpr2}) holds true, we exploit again
the universality of $\pi$ (but now the \"uniqueness" part).
We want to show that in the right hand special triangle of the diagram below, $i$ is the identity:
\[
\xymatrix@R=1pc@C=1pc{
&&& f^*P\ar@{-->}[rr]^{f}\ar@/^1pc/[dddddl]|(.4){f^*i} & & P  \ar@/^1pc/[dddddl]^i  \\
f^*P\ar@{-->}@/^0.8pc/[urrr]^\id\ar@{-->}[rr]^\Delta \ar[dd]_\pi && 
W \ar@{-->}[rr]\ar@{-->}[ur]\ar[dd]_{h^*\pi} && f^*P \ar[dd]_\pi\ar@{-->}[ur]_f  \\ 
&&&& & & & \\
X\ar@{-->}[rr]^\Delta\ar@/_1pc/[rrdddd]_\id && h^*X \ar@{-->}[rr]^h\ar[dd]_{f^*u} && X \ar[dd]_u & & & \\ \\
&& f^*P \ar@{-->}[rr]^f\ar[dd]_\pi && P &&&    \\ \\  
&& X }
\]
By the Beck-Chevalley and Frobenius laws, the left hand triangle with left side $\id_X\pi$ 
and right side $\pi(f^*i)$ is also special.
Since its top side is the identity, we get $\pi(f^*i) = \pi$. So $i$ and $\id_P$
both give a factorization of $\pi$ through $\pi$, and we get $i = \id_P$ as desired.

Lastly, we prove that (R) implies (AP). 
Suppose then that $u:X\to P$ is a stably opcartesian arrow over $f$.
We show that $u$ is part of an algebraic product $(\pi,u)$. 
The map $\pi:f^*P\to X$ is obtained by exploiting the universal property of $u$
(in fact, the \"existence" part).
It is the (unique) map such that $\pi(f^*u) = \Delta_!\id$, that is 
such that the triangle in the diagram below is special:
\[
\xymatrix@R=4pc@C=1.5pc{
X \ar@/_1pc/[ddr]_\id \ar@{-->}[rr]^\Delta  & & h^*X \ar[d]^{f^*u}\ar@{-->}[rr]^h && X \ar[d]^u \\ 
&  & f^*P \ar@/^.8pc/@{..>}[dl]^\pi \ar@{-->}[rr]_f & & P \\
& X &      }    
\]
To show that also the first of equations (\ref{algpr2}) holds true, we exploit again
the universality of $u$ (but now the \"uniqueness" part).
We want to show that in the bottom special triangle of the diagram below, $i$ is the identity:
\[
\xymatrix@R=1pc@C=1.5pc{
X \ar@{-->}[rr]^\Delta \ar@/_1pc/[ddddrr]_\id \ar@/^2pc/@{-->}[rrrr]^\id 
&& h^*X \ar[dd]_{f^*u}\ar@{-->}[rr]^h  & & X \ar[dd]^u \\ 
& \# & &@ & \\
&& f^*P \ar[dd]_\pi\ar@{-->}[rr]^f  & & P \ar@/^1pc/[ddddl]^i \\ \\
&& X \ar@/_/[ddr]_u & \# & & \\ \\
&& & P &      } 
 \]
The outer triangle is special (by pasting and compositions) and has an identity top side,
so that $iu = u = \id_Pu$. Thus, $i = \id_P$ as desired.
\epf

\begin{remark}
The present proof of theorem \ref{maincart} is the same as that in \cite{pisani2}
and follows, in the setting of unbiased symmetric multicategories, 
the proof given in \cite{pisani} for the classical case.
\end{remark}


\subsection{Cartesian cocartesian multicategories}
\label{Cartco}
Recall (see Section \ref{Coca}) that for any category $\C$, we have the associated cocartesian (or \"sequential")
symmetric multicategory $M:\dC\t\to\dPb$.
As sketched in example \ref{cartex2}, if $\C$ is enriched in commutative monoids then $\dC\t$ has a cartesian structure.
Conversely, it easy to see that a cartesian structure on $\dC\t$ gives an enrichment of $\C$ in $\cMon$.
In fact we have (\cite{pisani}):
\begin{proposition}
Cartesian structures on $\dC\t$ correspond to enrichments of $\,\C$ in commutative monoids.
\end{proposition} \epf
To state it in a more colorful way: 
\q\"cartesian + cocartesian = $\cMon$-enriched".

As a corollary, we find again the following well-known fact:
\begin{proposition}
A category $\C$, enriched in commutative monoids, has finite sums if and only if 
it has finite products and, if this is the case, they coincide and are absolute.
\end{proposition}
\pf
Since $\dC\t$ is cocartesian, it is representable if and only if $\C$ has finite sums.
Since $\dC\t$ is cartesian, by Theorem \ref{maincart} it is representable if and only it has universal products,
if and only if $\C$ has finite products; furthermore, the products are algebraic, hence preserved by any 
functor in $\cMon\-\Cat$.
\epf


\subsection{Comparison between cartesian and cocartesian multicategories}
\label{Comp}
It is interesting to compare cartesian and cocartesian multicategories: 
as the terminology suggests, they are ideally dual notions, though it is clearly a non-perfect duality.
Here are some similarities and differences.
\begin{enumerate}
\item
\begin{itemize}
\item
If a cartesian $M$ is representable, it is represented by universal products in $M$, which are also
products in the underlying category.

Representability is an opfibration condition for loose arrows, universal products are an opfibration condition
for what corresponds conceptually to loose arrows in the opposite direction.
\item
If a cocartesian $M$ is representable, it is represented by sums in the underlying category.

Representability and sums are opfibration conditions for loose arrows, cocartesianess is a fibration condition
for loose arrows.
\end{itemize}
\item
\begin{itemize}
\item
To be cartesian is a {\em structure} on a multicategory $M$; in fact, $M$ can support non-isomorphic
cartesian structures.
For instance, two non-isomorphic rings $R$ and $R'$ with the same underlying multiplicative monoid $M$
give two non-isomorphic cartesian structures on the same (cocartesian) multicategory $M\t$.

In fact, cartesian multicategories are the algebras for a monad that is not (even laxly) idempotent.
Roughly, the monad adds the opposite of the tight arrows to the loose part of $M$.
It doesn't seem to be directly associated to a tensor product on multicategories.

\item
To be cocartesian is a {\em property} on a multicategory $M$; in fact, 
cocartesian muticategories can be characterized as those symmetric multicategories with a \"central monoid"
(see \cite{pisani}), that is each object supports a monoid structure which commutes with all the arrows of $M$.

In fact, cocartesian multicategories are the algebras for an idempotent monad.
Roughly, the monad adds the tight arrows to the loose part of $M$ (and then will form the above mentioned central monoid).
It can be directly associated to the Boardman-Vogt tensor product on symmetric multicategories (see \cite{pisani}).

\end{itemize}

\end{enumerate}


\subsection{Cartesian multicategories as theories}
\label{Carth}
As mentioned in sections \ref{Unmon} and  \ref{Smc},
if $\dM \in \DP$ is the double prop corresponding to a symmetric multicategory $M:\dM \to \dPb$, 
the category of $\dM$-algebras is the fibered form of the usual category of $M$-algebras,
that is, of functors from $M$ to the symmetric multicategory of sets and several variable mappings.
Thus, a multicategory $M:\dM \to \dPb$ can be seen as a sort of theory, 
whose \"strict models" $M\to\Set$ can be captured in the fibrational form of
$\dM$-multicategories $M':\dM' \to \dM$ with a discrete opfibration condition on loose arrows.

Similarly, a {\em cartesian} multicategory $M:\dM \to \dPb$ can be seen as a sort of theory, 
whose \"strict models" $M\to\Set$ in $\sMlt^\cart$ can also be captured in the fibrational form of
$\dM$-multicategories $M':\dM' \to \dM$ with a discrete opfibration condition on loose arrows.

In this case, the multicategory $M\circ M':\dM' \to \dPb$ in $\sMlt$, 
also has a cartesian structure and the morphism $M':M\circ M'\to M$ in $\sMlt$ preserves the cartesian structures,
that is, it is a morphism in $\sMlt^\cart$.
For instance, if $M = R\t$ is the cartesian cocartesian multicategory associated to a rig $R$,
a model $M\to\Set$ is, as already remarked, an $R$-module $A$.
The corresponding multicategory $M\circ M'$ has the elements of $A$ as objects and 
linear combinations as arrows. Its cartesian structure encompasses the properties 
$ra + sa = (r+s)a$ and $0 = 0a$.


\subsection{Free cartesian multicategories}
\label{Further}
Since unbiased  cartesian muticategories are the algebras for the monad $(-)^\cart$ on $\sMlt$,
for any symmetric multicategory $M:\dM\to\dPb$, we have the free cartesian muticategory 
$M^\cart:\dEspan(\dM)\to\dPb$.
Its models $M^\cart \to \Set$ in $\sMlt^\cart$ correspond to the models $M\to \Set$ in $\sMlt$.
Here are some examples.

If $M:\dC\t\to\dPb$ is familial on $\C$, then $M^\cart$ is again familial:
\eq
M^\cart : \dC'\t\to\dPb
\eeq
where $\C'$ is the free $\cMon$-enriched category on $\C$,
that is, $\C'(A,B)$ is the free commutative monoid on $\C(A,B)$.
Indeed, in this case the action tight arrows is free, since the loose arrows of $M^\cart$ are the enhanced spans themselves. 

In particular, for a monoid $M$ we get the free rig $M'$ on $M$ of linear combinations with natural coefficients.
An action of $M$ on a commutative monoid $C$ extends uniquely to an $M'$-module on $C$.
In particular, for the terminal symmetric multicategory $1\t$ (which is the identity $\dPb\to\dPb$),
we get, as already noticed, the rig $\dN$ of natural numbers.

Recall that $U:\dBij\to\dPb$ is the terminal unary multicategory, whose algebras are sets (\ref{bijmon}).
Then $U^\cart$ is the cartesian theory for sets.
Explicitly, it can be seen as the submulticategory of $\dN\t$ 
whose arrows are families $I\to\{0,1\}$ with {\em exactly one} $1$.
The multiplication-composition and the sum-reindexing are as those of $\dN\t$.
Thus, in the passage from $U$ to $U^\cart$, one adds operations which become projections in the algebras: 
projections are the \"linear combinations" with the above mentioned families as coefficients. 
If we pass to the corresponding Lawvere theory, as for $\dN\t$ (or, more generally, $\dC\t$), we pass
from families to matrices.
What we get in this case are \"functional matrices" or \"functional relations", that is the usual Lawvere theory
$\Set\op$ for sets.
Note that $U^\cart$ is equivalent to $\dSec\to\dPb$, the multicategory associated to
the species of elements (\ref{ex4}).

Examples of non-free cartesian multicategories are those categories enriched in non-free commutative monoids.
For instance, the terminal symmetric multicategory $1\t : \dPb \to \dPb$ has clearly just one cartesian structure,
whose models are singletons (terminal sets).
Indeed, these are the modules for the terminal rig, in which $0 = 1$, so that $a = 1a = 0a = 0$ for any $a$.


\begin{refs}

\bibitem[Beilinson \& Drinfeld, 2004]{beilinson} A. Beilinson and V. Drinfeld (2004), {\em Chiral algebras}, 
American Mathematical Society Colloquium Publications, 51.

\bibitem[Benabou, 1985]{benabou2} J. Benabou (1985), Fibered categories and the foundations of naive category theory, 
{\em The Journal of Symbolic Logic} {\bf 50}, 10-37.

\bibitem[Bergeron et al. 1998]{bergeron} F. Bergeron, G. Labelle, P. Leroux (1998), 
{\em Combinatorial species and tree-like structures}, Encyclopedia of Mathematics and its Applications, vol. 67, 
Cambridge University Press.

\bibitem[Fiore et al., 2018]{fiore} M. Fiore, N. Gambino, M. Hyland, G Winskel (2018), 
Relative pseudomonads, Kleisli bicategories, and substitution monoidal structures
{\em Selecta Mathematica}, Springer.

\bibitem[Gambino  \& Kock, 1013]{gambino} N. Gambino and J. Kock (2013), Polynomial functors and polynomial monads,
{\em Math. Proc. Cambridge Phil. Soc.} {\bf 154}, 153-192.

\bibitem[Hermida, 2000]{hermida} C. Hermida (2000), Representable multicategories, 
{\em Advances in Math.} {\bf 151}, 164-225.

\bibitem[Janelidze  \& Street, 2017]{jane} G. Janelidze and R. Street (2017), Real sets, 
{\em Tbilisi mathematical journal} {\bf 10}, 23-49.

\bibitem[Joyal, 1981]{joyal} A. Joyal (1981), Une théorie combinatoire des séries formelles, 
{\em Adv. in Math.} {\bf 42}, 1-82.

\bibitem[Kelly, 2005]{kelly} G.M. Kelly (2005), On the operads of J.P. May, 
{\em Reprints in Theory and Appl. of Cat.} {\bf 13}, 1-13.

\bibitem[Lambert, 2021]{lambert} M. Lambert (2021), Discrete double fibrations, 
{\em Theory and Appl. of Cat.}  {\bf 37}, 671-708.

\bibitem[Lambert \& Patterson]{patterson} M. Lambert and E. Patterson (2024),
Cartesian double theories: a double-categorical framework for categorical doctrines, 
{\em Adv. in Math.} {\bf 444}, 109630.

\bibitem[Leinster, 2003]{leinster} T. Leinster (2003), {\em Higher operads, higher categories}, Cambridge University Press.

\bibitem[Lurie, 2017]{lurie} J. Lurie (2017), {\em Higher algebra}, available on the web.

\bibitem[May \& Thomason, 1978]{may} J. P. May and R. Thomason (1978), The uniqueness of infinite loop space machines,
 {\em Topology} {\bf 17}, 215-224.

\bibitem[Pisani, 2014]{pisani} C. Pisani (2014), Sequential multicategories, 
{\em Theory and Appl. of Cat.}  {\bf 29}, 496-541.

\bibitem[Pisani, 2022a]{pisani2} C. Pisani (2022a), Fibered multicategory theory, 
preprint available on arXiv.org.

\bibitem[Pisani, 2022b]{pisani3} C. Pisani (2022b), Operads as double functors, 
preprint available on arXiv.org.

\bibitem[Streicher, 2018]{benabou} T. Streicher (2018), 
Fibered Categories a la Jean  Benabou, preprint available on arXiv.org.

\end{refs}

\end{document}